\documentclass[11pt]{epiarticle}
\usepackage{epimath}

\setpapertype{A4}

\usepackage{hyperref}
\hypersetup{linkcolor=blue,citecolor=red,filecolor=green,urlcolor=magenta} 
\usepackage[english,french]{babel}
\addto\captionsfrench{%
   \def\figurename{Figure}%
   \def\tablename{Tableau}%
}
\usepackage[T1]{fontenc}
\usepackage[utf8]{inputenc}

\usepackage[pdftex]{graphicx}  % A utiliser avec pdflatex
\usepackage{caption}
\usepackage{apacite}

\usepackage{tabularx}
\newcolumntype{C}{>{\centering}X}

\usepackage[titletoc]{appendix}

\usepackage{graphicx}
\usepackage{pgf,tikz}%pour dessins en latex
\usepackage{pgfplots}
\pgfplotsset{compat=1.15}
\usepackage{mathrsfs}
\usepackage{array} 
\usetikzlibrary{arrows}

\newcommand{\R}{{\mathbb R}}

\title{Des croquis comme support de raisonnement et de changement de registres}
\titlemark{\'Etude du r\^ole des croquis et schémas}
\author{D. Grenier, C. Menini,  P. Sénéchaud \&\  F. Vandebrouck  avec la CIIU}
\authormark{CIIU }
\date{} % Empty date or tweak it according to your needs
\journal{\'EpiDEMES  -- (aaaa), ---} % Epijournal name
\acceptation{soumis le jj/mm/aaaayyyy, accepté le jj/mm/aaaa, version définitive le jj/mm/aaaa}
\begin{document}

\selectlanguage{french}

%%%%%%%%%%%%%%%%%%%%%%%%%%%%%%%
% Ajouter le titre du document
%%%%%%%%%%%%%%%%%%%%%%%%%%%%%%%
\maketitle

%%%%%%%%%%%%%%%%%%%%%
% Dédicace (si nécessaire)
%%%%%%%%%%%%%%%%%%%%%
\dedication{}

%%%%%%%%%%%%%%%%%
% Remerciements (si nécessaire)
%%%%%%%%%%%%%%%%%

\thanks{}

%%%%%%%%%%%%%%%%%%%%%%%%%%%%%%%%%%%%%%%%%%%%%%%%%%%%%%%%%%
% Ajouter le résumé, les mots-clés, la classification MSC (recommandée)
% Ne jamais ôter la section prelims, la conserver vide au besoin
%%%%%%%%%%%%%%%%%%%%%%%%%%%%%%%%%%%%%%%%%%%%%%%%%%%%%%%%%%
\begin{prelims}

\selectlanguage{english}

\def\abstractname{Abstract}
\abstract{We are convinced of the usefulness of sketches and diagrams during mathematical work but the observation is made in our practices that they are not spontaneously used by students.
In order to study the understanding and use of sketches by mathematics students, we designed and then proposed a test at different university levels. The test consists of five exercises. The first concerns different representation registers of a set of numbers, the second on a graphic proof of an implicative algebraic proposition and the last three on the graphic approch to the notions of injectivity, surjectivity, bijectivity in the context of the analysis. The sketches, proposed or requested in each exercise, are intended to be aids to changes of register and reasoning.
We present what motivated the choices and developments of the exercises then we analyze the results of these tests. In each case, we see difficulties in understanding and the sketches proposed, which leads to think that the sketch must be the subject of specific work at the university level.}


\medskip

\selectlanguage{french}
\def\abstractname{Résumé}
\abstract{Nous sommes convaincus de l’utilité des croquis et schémas lors du travail mathématique mais le constat est fait dans nos pratiques qu’ils ne sont pas spontanément  utilisés par les étudiants. 
Afin d’étudier la compréhension et l’utilisation de croquis par des étudiants en mathématiques, nous avons conçu puis proposé un test à différents niveaux universitaires. Le test comporte cinq exercices. Le premier porte sur différents registres de représentation d’un ensemble de nombres, le deuxième sur une preuve graphique d’une proposition algébrique implicative et  les trois derniers sur l'approche graphique des notions d’injectivité, surjectivité, bijectivité dans le cadre de l’analyse.  Les croquis, proposés ou demandés dans chaque exercice, se veulent être des aides aux changements de registre et au raisonnement. 
Nous présentons ce qui a motivé les choix et évolutions des exercices puis nous analysons les résultats de ces tests. On constate dans chaque cas des difficultés à la compréhension et à l'appropriation des croquis proposés, ce qui conduit à penser que le croquis doit être l’objet d’un travail spécifique au niveau universitaire.}

\def\keywordsname{Mots-Clés}
\keywords{Croquis, Schéma, Registre de représentation, Didactique des mathématiques}

\selectlanguage{french}

% Ajouter une table des matières (optionnelle)
\tableofcontents

\end{prelims}

\section{Introduction}

La Commission Inter IREM Université (CIIU) permet la rencontre d'universitaires en mathématiques et didactique des mathématiques qui partagent leurs diverses compétences sur l'enseignement. Un de ses thèmes de réflexion est l'étude du rôle des schémas et croquis qui sont des registres de représentation sémiotique importants pour la construction de connaissances
opératoires et pérennes. Notre réflexion, \`a son origine, concernait l'utilisation, dans nos pratiques, d'esquisses rapides comportant un codage de ce qu'on veut mettre en avant \`a un moment pr\'ecis d'un cours ou de le la résolution d'un exercice. Qu'en retiennent les étudiants et s'en emparent-ils? 

Le constat est fait que les croquis, les dessins et les schémas ne sont pas naturellement mobilisés
par une majorité d'étudiants ou d'étudiantes pour aborder un problème ou représenter une notion mathématique.
Nous faisons l'hypothèse de leur utilité en ce qu'ils peuvent contribuer au développement d'une image mentale sur laquelle s'appuyer et qu'ils apportent  des outils qui rendent opérationnelle la notion \`a app\'ehender. Ils peuvent aussi permettre à l’enseignant  de déceler des lacunes, à l’étudiant de conjecturer, vérifier et corriger les arguments de son raisonnement. Mais ils peuvent aussi être source d'ambiguïté dans l'échange entre étudiants et enseignants.

L'usage des dessins et des représentations  a d\'ej\`a fait l'objet d'\'etudes : en particulier dans l'articulation entre alg\`ebre lin\'eaire et géométrie \cite{GUE1},  et dans la différence entre figure et dessin \cite{LC1} et \cite{PAR1} 
ou encore dans les  travaux de R.Duval \cite{DUV} o\`u le concept de figure est reli\'e \`a la notion d'objet math\'ematique et o\`u le dessin en est une repr\'esentation. 

Une première analyse issue du travail de la CIIU se trouve  dans \cite{MS1} dans laquelle sont donnés des exemples o\`u les croquis, dessins et  autres repr\'esentations  peuvent aider les \'etudiants dans leurs apprentissages. 
Suite à ce premier travail, des expérimentations ont été menées  afin de tester la compréhension et le rôle de certains croquis dans la résolution de problèmes. 

Nous avons construit un test  que vous trouverez en annexe \ref{test}. 
Pour mettre en place ce test, nous avons travaillé en amont sur différentes versions,  ce qui a permis un  enrichissement de l'analyse {\it a priori }. 
%et l'élaboration d'une grille de lecture. La grille a été utilis\'ee par le groupes liaison Lycée  de l'IREM de Strasbourg pour l'analyse des copies locales.

Le test a été proposé au niveau L2\footnote{Deuxième année universitaire}  mathématiques à Strasbourg (19 copies), en L3\footnote{Troisième année universitaire} mathématiques \`a  Bordeaux (11 copies), en  L3 préprofessonalisation MEEF\footnote{Métiers de l'Enseignement de l'Education et de la Formation} à Limoges (7 copies) et en M1 MEEF\footnote{Premi\`ere ann\'ee de Master d'enseignement} à Poitiers (8 copies).
\`A chaque fois, le test a été posé sans préparation préalable afin de mesurer au plus près ce qui fait partie des connaissances disponibles des étudiants et de leurs habitudes de travail.

Nous présentons une analyse de ce test et les résultats des travaux menés par la CIIU. 

Soulignons ici l'importance des croquis dans d'autres sciences comme par exemple en sciences de la vie et en géographie\footnote{\url{https://actu-hg.hatier.fr/sites/default/files/Le\%20croquis\%20en\%20Geo_0.pdf}}.  
En scicence de la vie, le {\bf dessin} est une repr\'esentation la plus fid\`ele possible de l'objet \'etudi\'e avec tous ses d\'etails, un {\bf croquis} est un dessin simplifi\'e sans d\'etails  pour comprendre la constitution et les \'eléments \`a \'etudier et enfin {\bf  le schéma}  form\'e de formes et de fl\`eches  pour expliquer la structure ou le fonctionnement de ce que l'on repr\'esente\footnote{\url{https://louisa-paulin.ecollege.haute-garonne.fr/espaces-pedagogiques/svt/fiches-methodologiques/comparaison-dessin-croquis-schema-46406.htm}}.
%{\color {purple} A regarder aussi 
% \url{https://www.alloprof.qc.ca/fr/eleves/bv/sciences/les-dessins-et-schemas-scientifiques-s1517} }
% {\color {purple}  En géographie.... } }. 
%En géographie aussi, le croquis (ou carte) a une place importante\footnote{\url{https://actu-hg.hatier.fr/sites/default/files/Le\%20croquis\%20en\%20Geo_0.pdf}}. 
 
 %Précision du vocabulaire :
%\begin{itemize}
%\item
%dessin : représentation soignée,
%\item
%croquis : représente une idée,
%\item
%schéma : les relations qui sont importantes pour l'étude sont matérialisées,
%\item
%diagrammes : {\color{red} on n'en a pas vraiment parlé ... on met tout de même ?}
%\end{itemize}

\begin{figure}[!ht]
%\centering
\includegraphics[width=8cm]{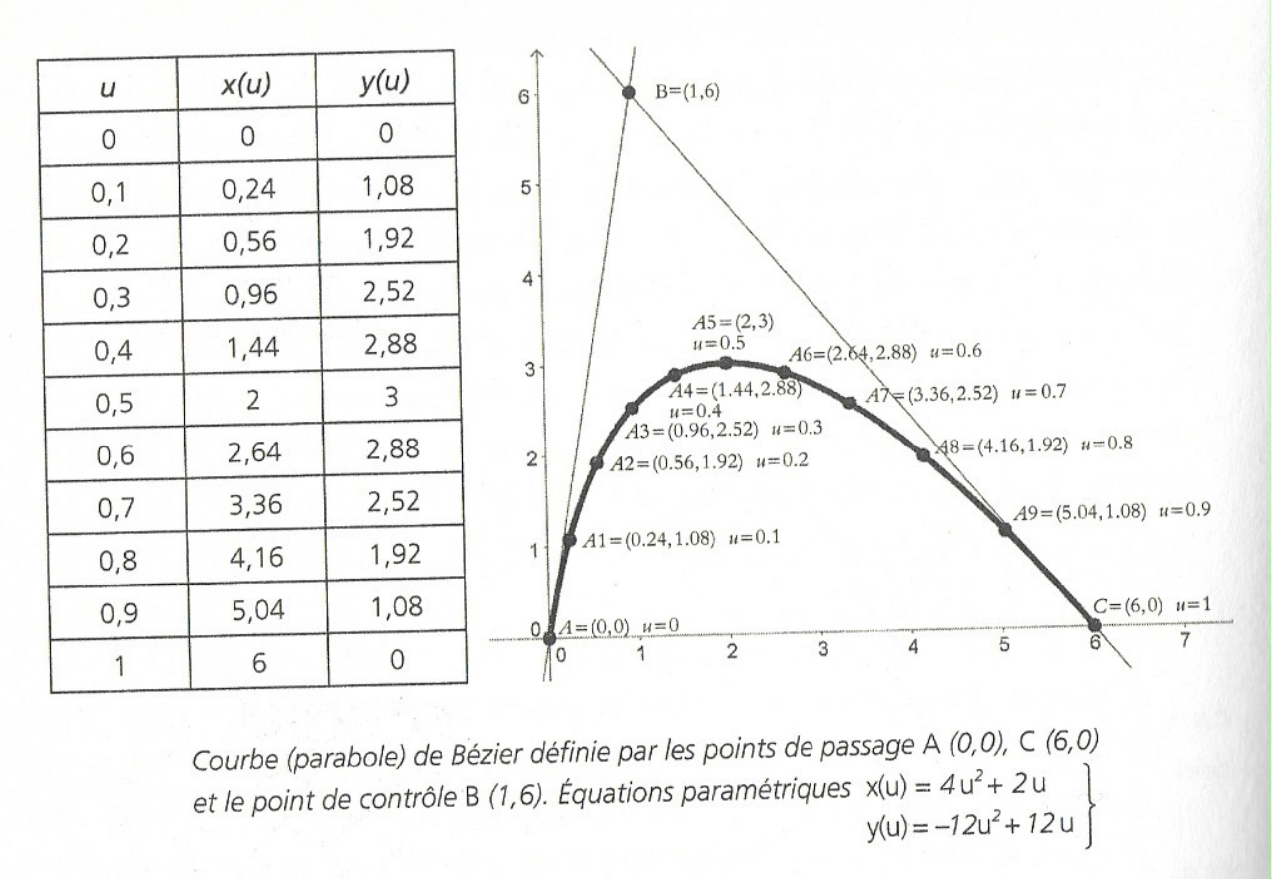} \includegraphics[width=8cm]{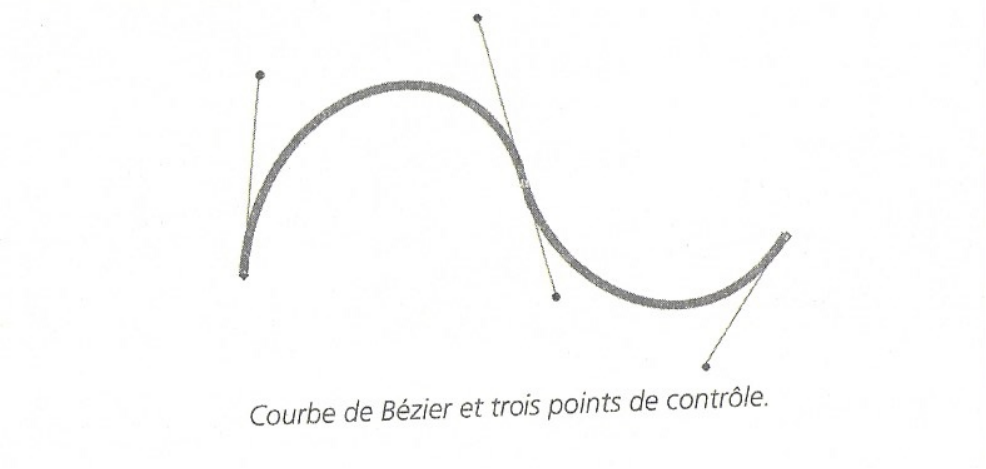} 
\caption{Un dessin / Un croquis  }\label{fig:bez1}
\includegraphics[width=8cm]{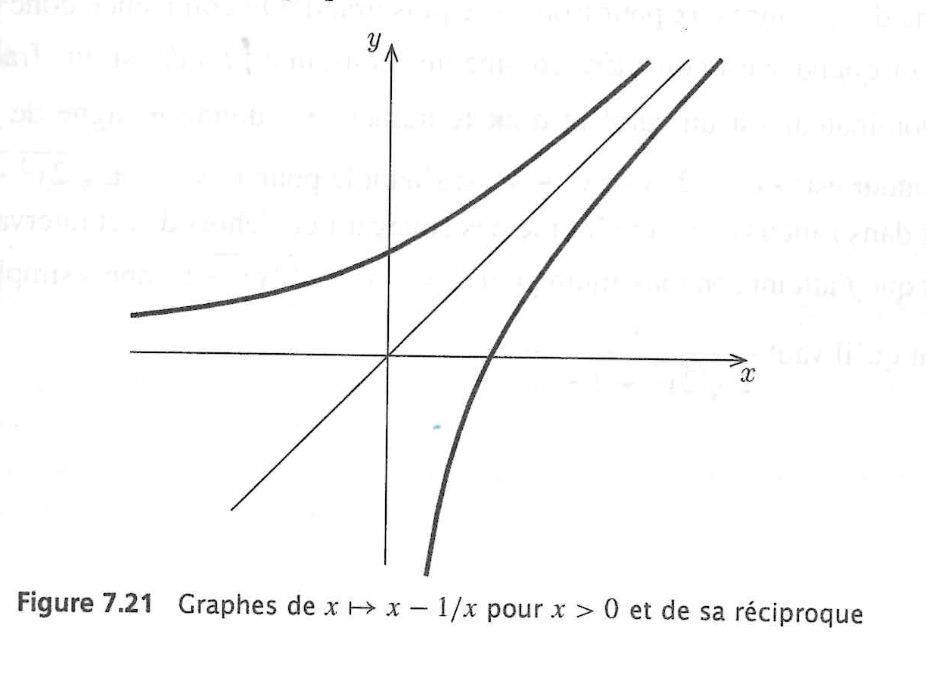} \includegraphics[width=8cm]{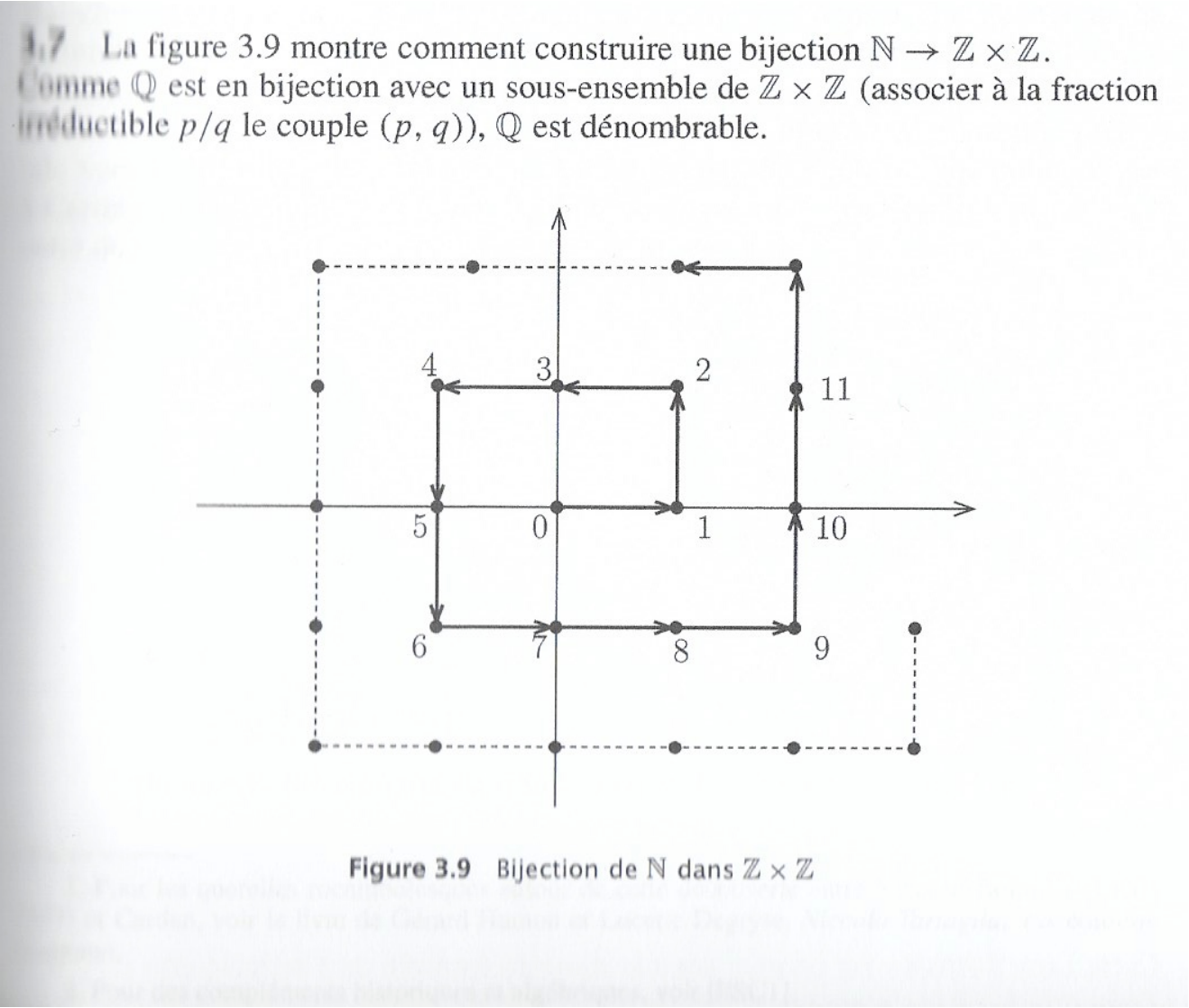} 
\caption{ Des croquis  }\label{fig:sansleg1}
\end{figure}

En math\'ematiques,  on peut reprendre la classification des sciences et vie de la terre~: un dessin est alors une repr\'esentation la plus compl\`ete possible de l'objet, un croquis est une représentation qui cherche \`a illustrer une  caract\'eristique ou propri\'et\'e particuli\`ere de l'objet d'\'etude et un schéma s'attache \`a  montrer  des relations, ou des articulations importantes pour l'\'etude de l'objet considéré (voir figure \ref{fig:bez1} et figure \ref{fig:sansleg1}).\ On pourra \'egalement noter que le sch\'ema permet de mettre en avant des propri\'et\'es g\'en\'eriques de l'objet.

En plus de savoir manipuler de telles représentations, il est également essentiel de savoir les articuler et en changer si besoin pour trouver celle la mieux adaptée à ce que l'on étudie ou à ce que l'on veut expliquer. Savoir changer de registres est un élément important dans l'apprentissage et la transmission des mathématiques : le test présenté permet d'étudier en partie ce savoir-faire. 
 
Ci-dessous, on peut cataloguer la repr\'esentation de gauche dans la figure \ref{fig:bez1}  comme un dessin, qui est assez pr\'ecis et contient les \'el\'ements de construction. \`A droite, il s'agit d'un croquis qui consiste \`a montrer les changements de concavit\'e et les points de contr\^ole\cite{COURB}.

Dans la figure \ref{fig:sansleg1}, nous avons deux croquis : celui de gauche  illustrant la propri\'et\'e de symétrie entre le trac\'e du graphe d'une fonction et de celui de sa r\'eciproque dans un rep\`ere orthonorm\'e, celui de droite illustrant la construction d'une bijection entre deux ensembles d\'enombrables\cite{ESC}.\\
%( Escoffier  Toute l'Analyse de Licence  Dunod ou encore dans  Courbes p\'erilleuses le monde est math\'ematiques) 
%Dire ici que pour avoir confiance dans un croquis il faut déjà avoir fait plusieurs dessins plus soignés, que le schéma ne peut être fait que si on a déjà compris ce qui est important d'un point de vue mathématique ?}
%{\color {purple} vu les d\'efinitions donn\'ees en SVT la finalit\'e d'un croquis n'est pas la m\^eme qu'un dessin mais clairement il s'agit de faire d'abord des dessins pour pouvoir faire des croquis car il faut avoir compris la constitution des \'el\'ements. } 

\begin{figure}[!ht]
\centering
\includegraphics[width=5.5cm]{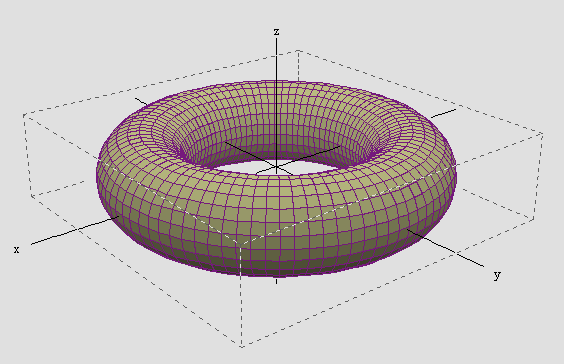} \includegraphics[width=7cm]{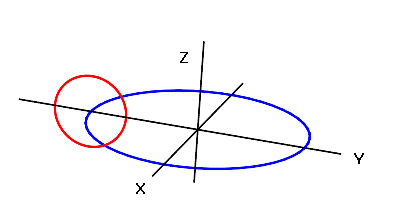} \includegraphics[width=3.5cm]{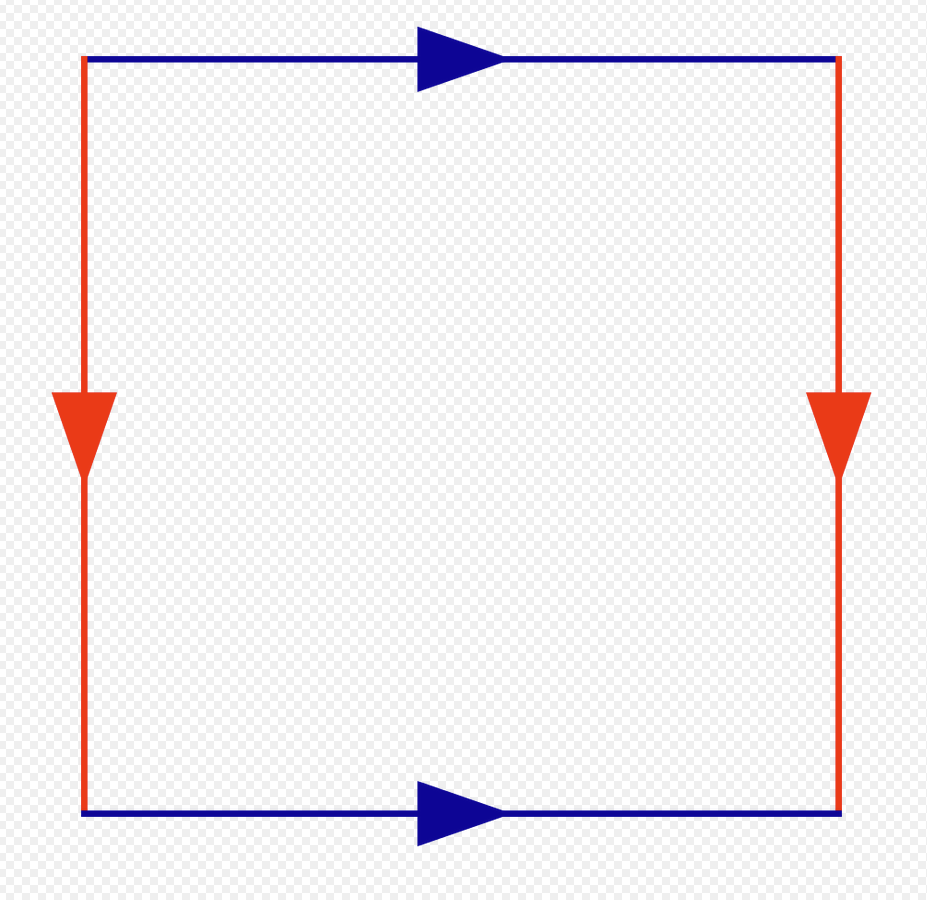}
\caption{Un dessin / Deux schémas    }\label{fig:Tore}
\end{figure}

\pagebreak
Dans la figure \ref{fig:Tore}, nous proposons trois repr\'esentations du tore : un dessin et deux sch\'emas. 
Celle de gauche est le dessin d'un tore (Wikipedia: domaine public),  celle du milieu sch\'ematise la construction par rotation du tore et celle de droite  sch\'ematise le 2-Tore par recollement des c\^ot\'es oppos\'es d'un carr\'e (Wikipedia: domaine public).

%%%%%La présentation ci-dessous ne compile pas !!!!%%%%
%\begin{figure}[!ht]
%\centering
%%\begin{subfigure}[a]{0.3\textwidth}
%\begin{subfigure}[]{}
%         %\centering
%         \includegraphics[width=5.5cm]{Tore2}
%         \caption{Représentation de la surface}
%         \label{fig:toredessin}
%     \end{subfigure}
%     \hfill
%     \begin{subfigure}[]{}
%         %\centering
%        \includegraphics[width=7cm]{Tore0}
%         \caption{Mise en évidence de la construction par rotation}
%         \label{fig:torerotation}
%     \end{subfigure}
%     \hfill
%     \begin{subfigure}[]{}
%         %\centering
%         \includegraphics[width=3.5cm]{Tore1}
%         \caption{Le tore vu comme $\R/\Z\times\R/\Z$}
%         \label{fig:torequotient}
%     \end{subfigure}
%\caption{Sch\'ema /dessin   }\label{fig:Tore}
%\end{figure}
% 

%Par Cham sur Wikipédia français — Transféré de fr.wikipedia à Commons par Bloody-libu utilisant CommonsHelper., Domaine public, https://commons.wikimedia.org/w/index.php?curid=18041326

%Par Theon — travail personnel (own work), inspiré de Image:Klein_Bottle_Folding_1.svg, CC BY-SA 3.0, https://commons.wikimedia.org/w/index.php?curid=4278162

Nous pouvons faire le parall\`ele avec les repr\'esentations  mobilis\'ees aupr\`es des \'el\`eves et des \'etudiants  lors de l'\'etude d'une fonction :
on peut leur demander une repr\'esentation tr\`es soign\'ee (en utilisant \'eventuellement un logiciel), ou encore l'allure de la courbe avec quelques points significatifs et quelques tangentes ou encore le tableau des variations qui est une autre repr\'esentation.

Le test proposé est composé de deux parties, une première qui conduit à un travail sur les notions de valeur absolue, distance et intervalles et une deuxième  dans laquelle le travail demand\'e fait appel aux  notions d'injectivité et de surjectivité pour des fonctions num\'eriques \`a variable r\'eelle.
Ces notions ont été choisies en raison de leur place prépondérante au début d'un  cursus universitaire classique (Licence généraliste en sciences). Pour chacune de ces deux parties, constitu\'ee de plusieurs exercices,  nous allons motiver le choix des notions traitées, donner une analyse {\it a priori}  qui correspond à ce à quoi on pourrait s'attendre  tant en terme de difficultés que de réussites et proposer une analyse des travaux des étudiants. Dans la conclusion de ces études nous proposons également des pistes d'amélioration des tests en fonction de ce qui est constaté.

\section{Diff\'erents registres de repr\'esentation  d'un ensemble de nombres }

\subsection{Présentation de l'exercice 1 et analyse {\it a priori } }

Dans le  premier exercice proposé, il s'agit de remplir des cases dans un tableau où chaque colonne correspond \`a un registre  différent de repr\'esentation d'un ensemble de nombres. Il s'agit d'évaluer la capacité à passer d'un registre à un autre en mobilisant les notions de valeur absolue, d'intervalle et de distance entre deux réels.

La valeur absolue est  au programme au début de lycée général et technologique. On peut lire dans les programmes de seconde (annexe du BO Numéro 1 du 22/01/2019) partie \og nombres et calculs\fg{} : 
\begin{itemize}
\item
{\it  La notation de la valeur absolue est introduite pour exprimer la distance entre deux nombres
réels et caractériser les intervalles de centre donné.} (alinéa manipuler les nombres réels)
\item
{\it Notation $|a|$. Distance entre deux nombres réels. Représentation de l’intervalle $[a - r , a + r]$ puis caractérisation par la condition
$|x - a| \leq  r$} (alinéa  contenus).
\end{itemize} 
 
La mise en relation entre ces deux notions  (valeurs absolue et distance entre deux réels) est indispensable en L1\footnote{Première année d'université} dès lors que sont travaillées les notions de limite et de continuité. Nous émettons, par ailleurs, l'hypothèse que si ce lien est bien construit en première année, il permet d'aborder, plus tard, certaines notions de topologie avec plus de facilité. 

% \noindent La première ligne du tableau est la suivante \\ 
%\begin{figure}[!ht]
%\centering
 %\includegraphics[width=14cm]{entete-colonne} 
%\caption{les différents registres présentés}\label{ligne1}
%\end{figure} 

On identifie le point d'abscisse $x$ avec le r\'eel $x$. Les colonnes correspondent \`a  des registres diff\'erents. Les colonnes 1 et 2 font appel au registre alg\'ebrique, la colonne 3 au registre graphique et les colonnes 4 et 5 au registre g\'eom\'etrique, comme l'indique les intitulés des colonnes. 
Le tableau compte huit lignes à remplir  et sur chaque ligne, une case est pré-remplie : aux étudiants de compléter, dans l'ordre qu'ils préfèrent et en numérotant les cases dans l'ordre de remplissage.
 Sur la première ligne, la première colonne est pré-remplie (par $x \in ]6,15[$), sur la deuxième ligne c'est la case de la deuxième colonne  (avec $-\sqrt 2 \leq x \leq \sqrt 2$). Sur la troisième ligne, c'est la colonne \og croquis \fg$\;$ qui est pré-remplie. Les colonnes \og valeur absolue \fg$\;$  et \og distance \fg$\;$ sont pré-remplies sur les lignes  4 à 7 et sur la dernière ligne un croquis est donné. 
 
Les cases pre-remplies  fournissent donc des indices pour remplir la totalité du tableau et, sur chaque  ligne, une colonne est remplie de manière à avoir au moins une case par colonne, sur l'ensemble du tableau : 
  
{ \tiny  \begin{center}
\renewcommand{\arraystretch}{3}
\begin{tabular}{|p{1.5cm}|p{2cm}|p{3.5cm}|p{1.5cm}|p{1.5cm}|}

\hline 

\begin{minipage}{2.3cm}ensemble fini\\ de réels ou \\ intervalle(s) \end{minipage} & \begin{minipage}{3cm}\'egalit\'es ou \\ in\'egalit\'es \end{minipage} & \begin{minipage}{7cm}croquis ou sch\'ema de \\ l'ensemble des $x$ considérés \end{minipage} & \begin{minipage}{2.3cm}valeur \\ absolue \end{minipage} & distance  \\ 

\hline 

$x\in]6,15[$ & &  &   & \\ 

\hline 

 & $-\sqrt{2}\leq x\leq \sqrt{2}$ &  &   & \\ 
 
 \hline

 &  & 
\begin{tikzpicture}[line cap=round,x=1.5cm,y=1cm, scale=0.5]
%\clip(0,-3) rectangle (6,1);
\draw[dashed,color=gray] (1,0)--(5,0);
\draw[line width=1.5,color=blue] ((1.5,0)--(4.5,-0);
\draw[line width=1.5,color=blue] (1.45,0.1)--(1.5,0.1)--(1.5,-0.1) node[below] {$-11$}--(1.45,-0.1);
\draw[line width=1.5,color=blue] ((4.55,0.1) --(4.5,0.1)--(4.5,-0.1) node[below] {$0$}--(4.55,-0.1);
\draw[line width=1,color=gray] (3,0.1)--(3,-0.1) node[below] {\color{blue}$-\frac{11}{2}$};

\draw[color=red,-triangle 60] ((3,.2) to[out =150, in=60] node[midway,above] {$-\frac{11}{2}$} (1.5,.2);
\draw[color=red,-triangle 60] ((3,.2) to[out =30, in=120] node[midway,above] {$+\frac{11}{2}$} (4.5,.2);

\end{tikzpicture}
  &  &  \\ 

\hline 

& & & & $d(x;\frac{3}{2})=1$\\

\hline

& & & $|x+5|=\pi$ & \\

 \hline 

&  &  & & $d(x;-4)\leq 5$ \\ 

\hline

 &  &  & $|x-3|\leq \dfrac{1}{2}$  & \\ 
 
\hline

& & 
\begin{tikzpicture}[line cap=round,x=1.5cm,y=1cm, scale=0.5]
%\clip(0,-3) rectangle (6,1);
\draw[dashed,color=gray] (1,0)--(5,0);
\draw[line width=1.5,color=blue] ((1,0)--(2,-0);
\draw[line width=1.5,color=blue] ((4,0)--(5,-0);
\draw[line width=1.5,color=blue] (1.95,0.1)--(2,0.1)--(2,-0.1) node[below] {$-\frac{5}{3}$}--(1.95,-0.1);
\draw[line width=1.5,color=blue] ((4.05,0.1) --(4,0.1)--(4,-0.1) node[below] {$\frac{7}{3}$}--(4.05,-0.1);
\draw[line width=1,color=gray] (3,0.1)--(3,-0.1) node[below] {\color{blue}$\frac{1}{3}$};

\draw[color=red,-triangle 60] ((3,.2) to[out =150, in=60] node[midway,above] {$-2$} (2,.2);
\draw[color=red,-triangle 60] ((3,.2) to[out =30, in=120] node[midway,above] {$+2$} (4,.2);

\end{tikzpicture}
 & & \\

\hline 
\end{tabular} 
\end{center}

} 
 
Pour construire cet exercice nous avons choisi de ne pas pré-remplir une ligne complète pour que le remplissage ne se fasse pas par mimétisme, de ne pas mettre que des nombres entiers et des intervalles centrés en 0. Nous avons aussi fait appel à des ensembles réduits à deux éléments mais aussi à des réunions disjointes d’intervalles, en ce qui concerne les ensembles de solutions. Et enfin nous avons fait intervenir des inégalités larges et strictes.

Nous faisons l'hypoth\`ese que les croquis sur lesquels sont aussi représentées des données opératoires (point milieu, demi-longueur de l'intervalle)   favorisent la mise en relation entre les deux colonnes de gauche et les deux colonnes de droite. C'est pourquoi nous les avons positionnés dans la colonne centrale.  

Le test a été proposé en L2 mathématiques  à Strasbourg en deux versions : une avec et une  sans représentation graphique (respectivement 10 et 9 tests), et cela pour avoir des indicateurs sur l'impact de la présence des croquis sur les réponses. 

Nous pouvons remarquer que, sur une ligne donnée,  le remplissage,  de la colonne 1 à partir de la colonne 2  et réciproquement met en jeu une simple ré-écriture des informations.  Cette remarque est également valable  pour remplir la colonne 4 à partir de la 5 et réciproquement.
Si l'étudiant connaît  l'équivalence entre  l'appartenance à une partie de $\R$, et la notion d'encadrement alors il saura remplir les deux premières colonnes  et si il connaît le lien entre valeur absolue et distance, de même il saura partir des informations de la colonne 4 pour remplir la colonne 5 et réciproquement. On peut s'attendre à retrouver l'erreur classique $d(x,a) = \vert x+a \vert$.

Une case pré-remplie sur la colonne \og croquis \fg, en position médiane, devrait permettre de remplir les colonnes de droite, comme de gauche, puisque  des informations nécessaires y sont représentées. Il s'agit pour l'étudiant de  détecter les informations sur le croquis.\\
 Inversement pour remplir cette colonne  \og croquis \fg, le travail n'est pas le même selon le pré-remplissage donné :  si c'est à partir des informations contenues dans les colonnes de droite c'est juste traduire ces informations en terme de croquis, en revanche si c'est  à partir des informations contenues dans les colonnes à gauche, cela nécessite de trouver le point milieu de l'intervalle et sa demi-longueur. Transférer l'information de gauche à droite, dans le tableau, nécessite donc un changement de registre plus élaboré qu'une simple ré-écriture, pour lequel le passage par la représentation graphique peut donner des pistes. 
 Nous nous attendons à relever des erreurs entre les inégalités larges et strictes, et nous sommes conscients  que l'égalité peut être déroutante.
 
\subsection{Analyse {\it a posteriori }  }
\begin{figure}[!ht]
\centering
\includegraphics[width=7cm]{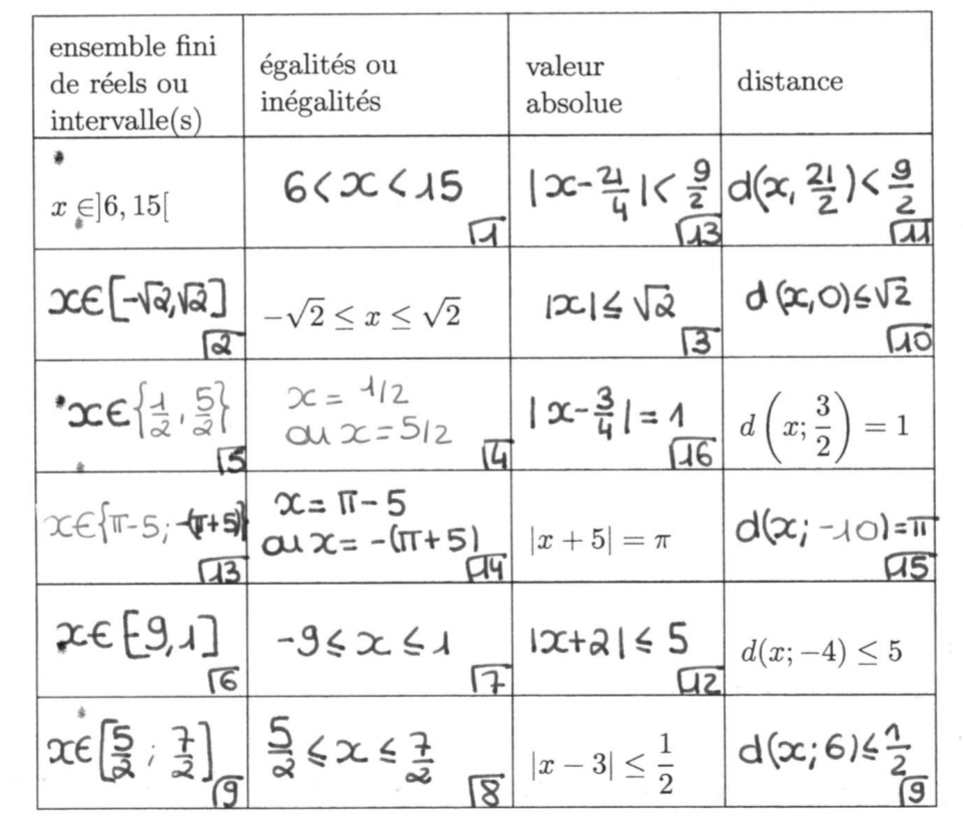} 
\caption{Souvenir \og qu'il faut prendre $1/2$\fg{} (L2-Strasbourg)}\label{fig:erreur1}
%\end{figure}
%\begin{figure}[!ht]
\centering
\includegraphics[width=7cm]{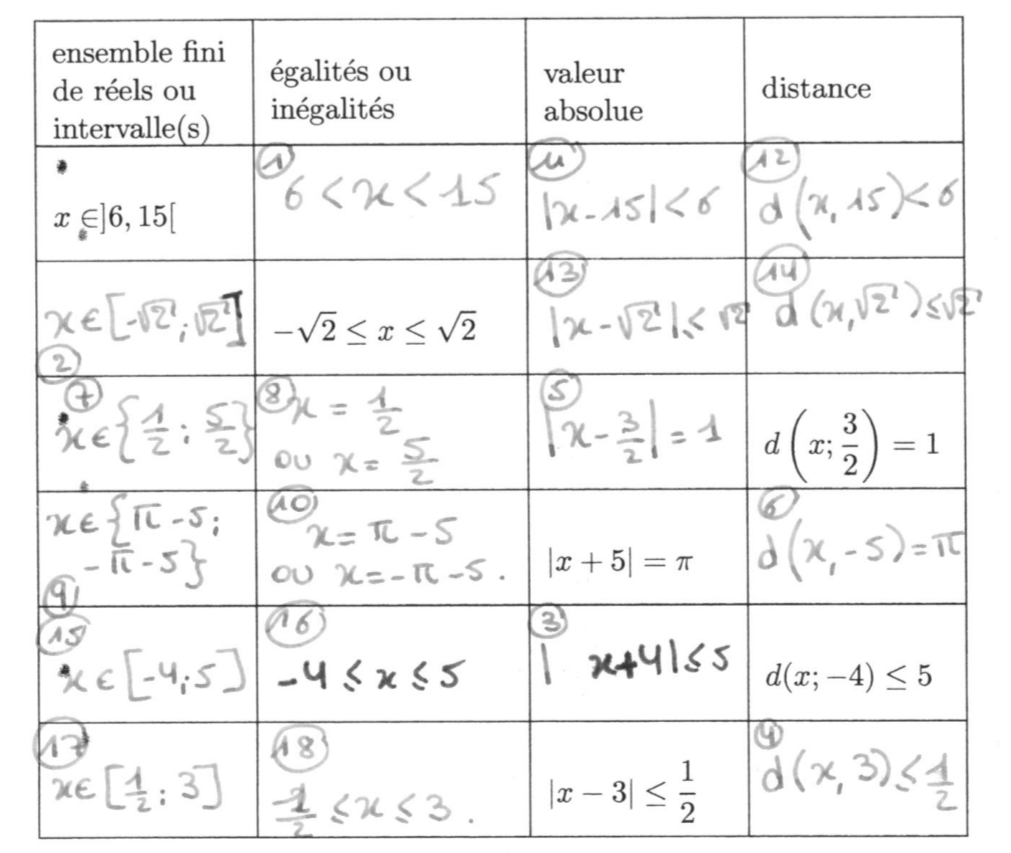} 
\caption{Un exemple où le lien n'est pas fait en abscence de croquis (L2-Strasbourg)}\label{fig:erreur2}
\end{figure}

Avant d'étudier l'influence de l'utilisation des graphiques, nous avons étudié les 9 copies de Strasbourg (L2)  où le tableau était donné sans la colonne  \og croquis \fg$\,$. Cette étude nous conforte  dans l'idée qu'elle pourrait être utile à certains étudiants. Dans les copies ci-dessous (figure \ref{fig:erreur1} et figure \ref{fig:erreur2})
les deux colonnes 1 et 2  d'une part et les colonnes 3 et 4 d'autre part, sont clairement remplies sans mise en relation entre celles de droite et celles de gauche.

\begin{figure}[!ht]
\begin{center}
\includegraphics[width=10cm]{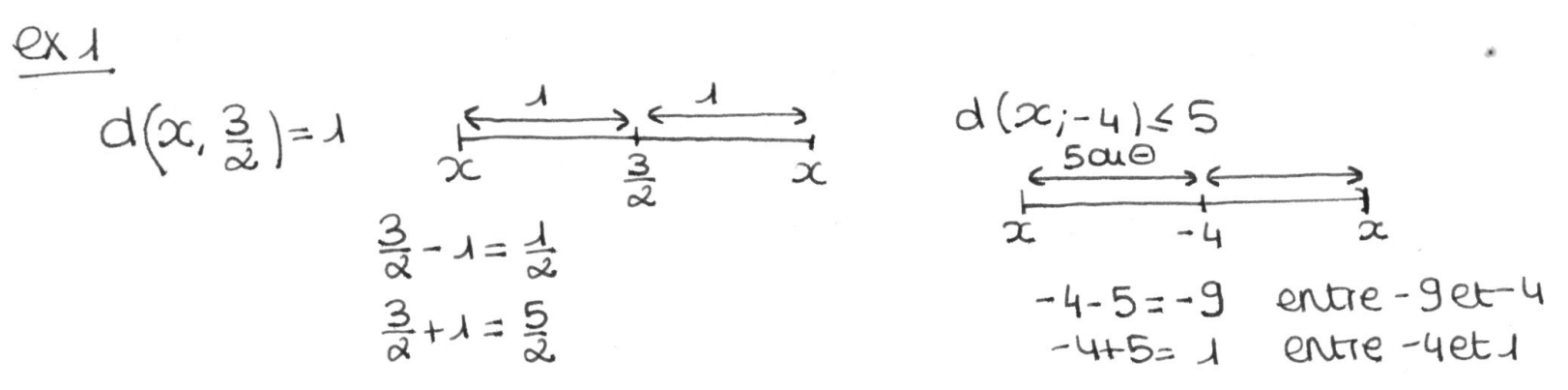} 
\caption{Schéma fait au brouillon pour un test sans croquis (L2-Strasbourg)}\label{fig:brouillon}
\end{center}
\end{figure}

Sur le brouillon d'une de ces  9 copies sans colonne  \og croquis \fg{} on constate le recours au schéma pour traiter le cas d'égalité (figure \ref{fig:brouillon}).

L'analyse  des 36  copies  (11 de Bordeaux, 7 de Limoges, 8 de Poitiers et 10 de Strasbourg) sur la version du test comportant une colonne \og croquis \fg$\,$ a été conduite de la façon suivante. Après analyse des deux premières colonnes et  pour les copies où les réponses sur ces colonnes sont correctes, nous avons étudié des points de vigilance  et nous nous sommes posés les questions suivantes. 
\begin{itemize}
\item
Le traitement de la colonne avec \og croquis \fg   : est-il correct ? erroné ? Quel type d'erreurs repère-t-on ? 
\item
L'utilisation du croquis  :  en
cas de réussite des croquis, les colonnes 4 et 5 sont-elles traitées correctement ?
\item
Le lien éventuel entre la réussite des colonnes de droite et la qualité du croquis : le croquis est-il détaillé ? contient-il des informations superflues ? 
\item
Le cas d'égalité est regardé à part.\\
\end{itemize}

En accord avec l'analyse {\it a priori}, on note une cohérence des réponses entre  les colonnes 1 et 2  : il s'agit d'une simple réécriture,  capacité mobilisable dans le registre algébrique.
La réécriture entre les colonnes 4 et 5 fonctionne également correctement. Toutefois, on retrouve l'erreur classique, signalée dans l'analyse {\it a priori} , $d(x,a)=|x+a|$.  Pour ceux qui ne savent pas passer de la distance à la valeur absolue,  le schéma n'est d'aucun secours. 

\begin{figure}[!ht]
\centering
\includegraphics[width=9cm]{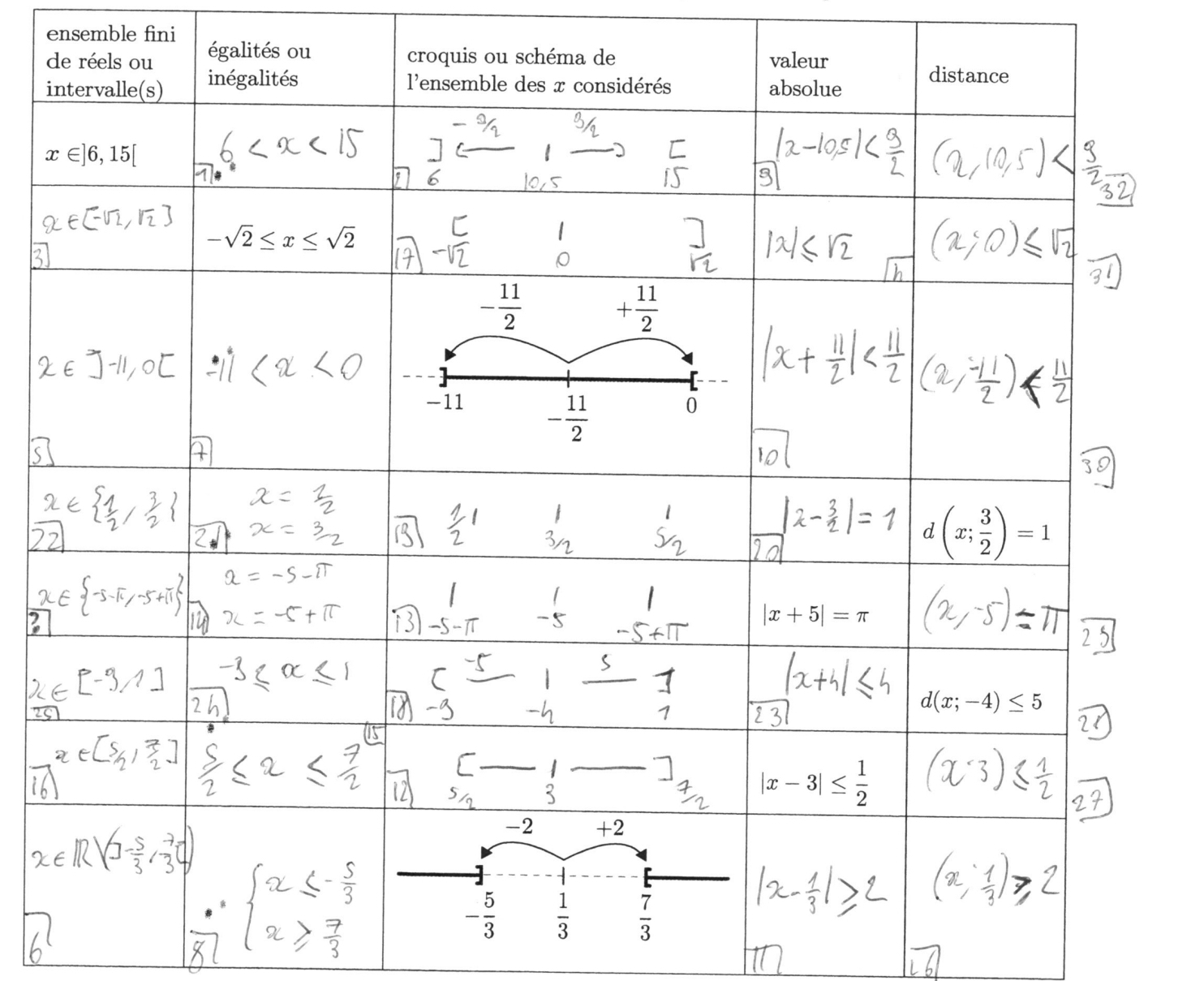} 
\caption{Exemple avec uniquement les données nécessaires (L2-Strasbourg)}\label{fig:sansaxe}
\end{figure}

Sur certaines copies correctement remplies,  le schéma est fait en dernier. Dans ce cas, le schéma ne semble  pas nécessaire à l'étudiant. La colonne \og croquis \fg,  qui correspond à une partie opérationnelle de l'activité, permet aux étudiants de construire et de représenter le milieu, mais cette compétence peut-être mobilisable sans cette colonne.
Par ailleurs, les étudiants qui réussissent le mieux font les tracés les plus approximatifs (certains ne dessinent même pas l’axe)  et mettent uniquement les données nécessaires, voir la figure \ref{fig:sansaxe}. On notera qu'il manque le $d$ de distance dans la dernière colonne.
\newpage 
Nous remarquons des confusions entre inégalités larges et strictes et la difficulté pour certains étudiants de représenter un ensemble à deux éléments.

\begin{figure}[!ht]
\begin{center}
\includegraphics[width=10cm]{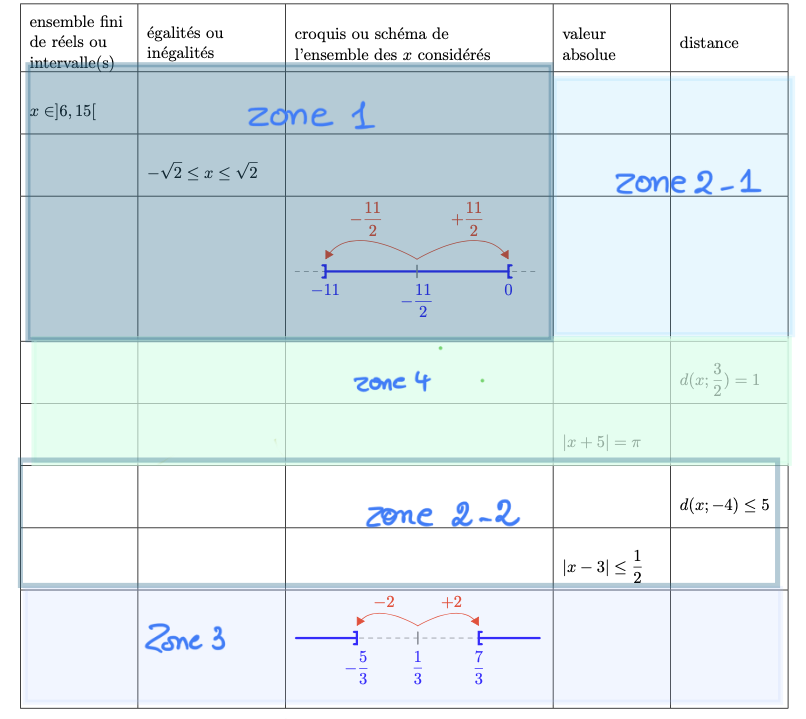} 
\caption{Zones de réussite}\label{fig:zones}
\end{center}
\end{figure}

En terme de réussite nous proposons une analyse plus détaillée qui partitionne le tableau en zones de réussite des étudiants  comme le montre la  figure \ref{fig:zones}. 

\begin{itemize}
\item
La zone 1 : les étudiants n'ont, en général, pas de difficultés à la remplir (seulement pour 2 copies sur les 36).   L'écriture de l'appartenance à un intervalle est bien comprise et correctement utilisée. Pour remplir la colonne centrale, celle du schéma, il s’agit de trouver le milieu de l'intervalle donné. Une fois cette information connue, les deux colonnes de droites sont faciles à remplir si les notions de valeur absolue et de distance sont connues. Cette partie est la zone 2-1.
\item
Les zones 2-1 et 2-2 :  elles ne sont pas toujours remplies en raison de la méconnaissance du lien entre valeur absolue et distance (dans 10 copies sur 36 toutes les cases de ces zones ne sont pas remplies). Sur ces lignes 6 et 7, le milieu de l'intervalle est donné et les bornes de l’intervalle sont données en terme de distance relativement au milieu. Le schéma en est une simple traduction. Lorsque le schéma est fait, les colonnes 1 et 2 sont remplies, sinon, non.  \`A partir de la distance ou de la valeur absolue (colonnes 4 et 5), les étudiants arrivent à remplir les colonnes 1 et 2. La colonne de la valeur absolue n'est pas toujours remplie, lorsque la distance est donn\'ee (figure~\ref{fig:zone22}).

\begin{figure}[!ht]
\begin{center}
\includegraphics[width=14cm]{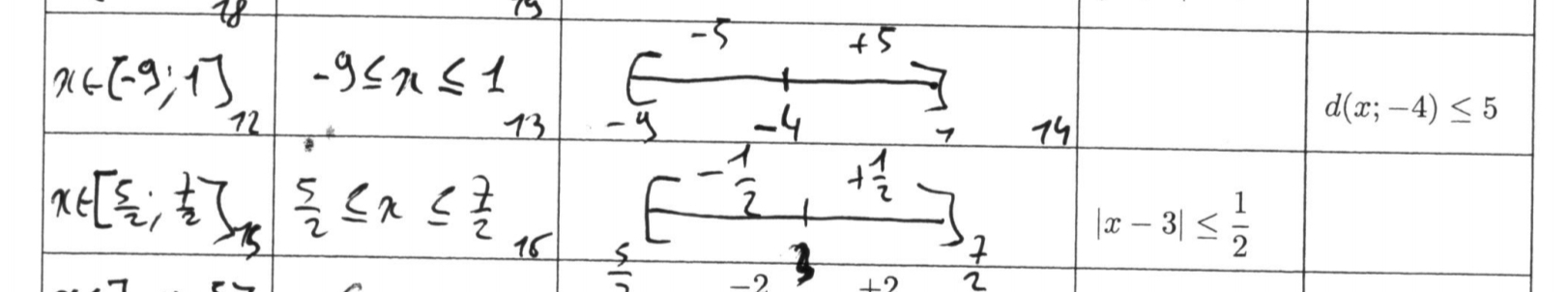} 

\vspace{0.5cm}

\includegraphics[width=14cm]{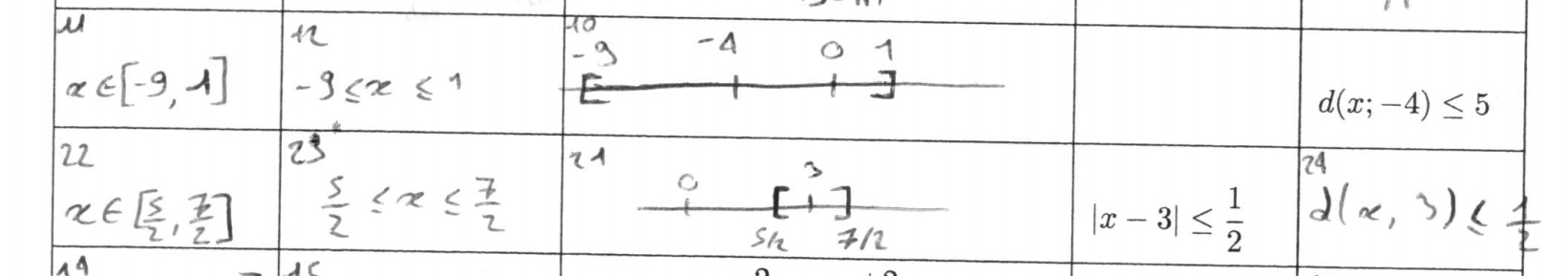} 
\caption{Zoom zone 2-2 : deux productions (L2-Strasbourg)}\label{fig:zone22}
\end{center}
\end{figure}

\vspace{5cm}

\item
Zone 3 : la dernière ligne est rarement complète en particulier sur les colonnes 4 et 5  (valeur absolue/ distance) en L2 surtout (4 sur les 10  de L2 Strasbourg remplissent correctement et 20  sur les 26  copies de L3 et Master1 ) : les étudiants arrivent à traduire le dessin en terme d’appartenance à une union d’intervalles ou à l’écrire avec des inégalités mais, sur les copies montrant le plus de fragilité, ni la colonne 4 ni la colonne 5 ne sont remplies (figure \ref{fig:zone3}).
\begin{figure}[!ht]
\begin{center}
\includegraphics[width=14cm]{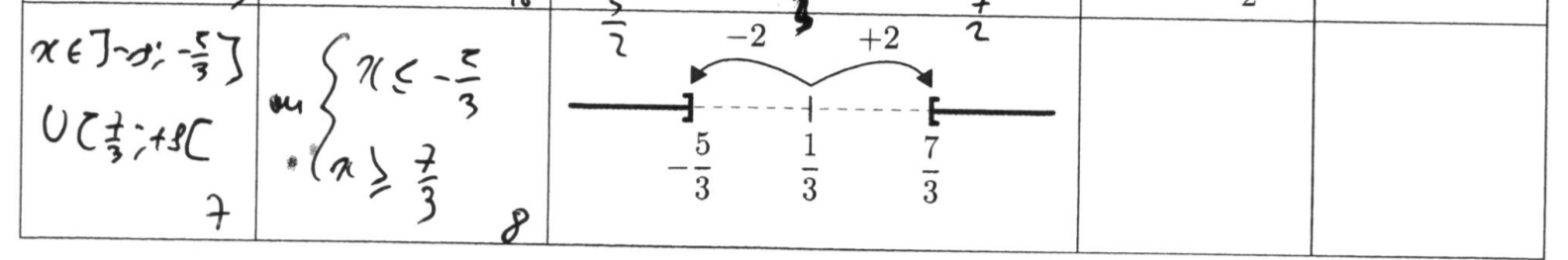} 

\vspace{0.5cm}

\includegraphics[width=14cm]{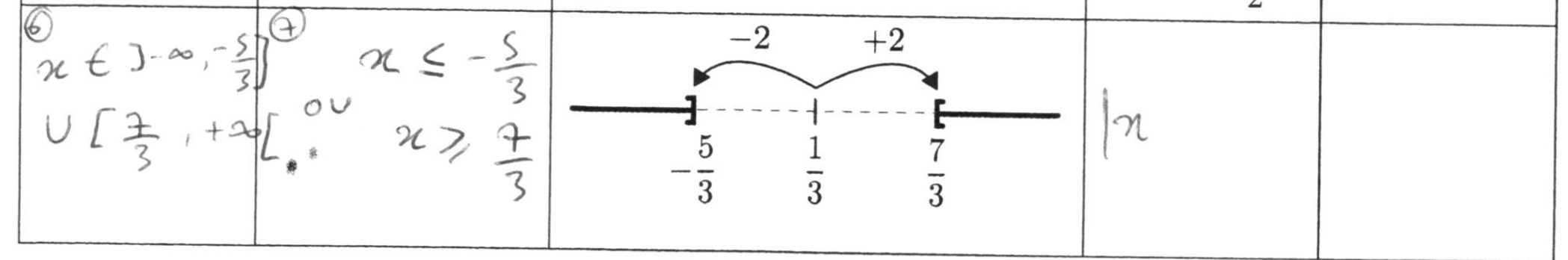} 
\caption{Zoom zone 3 : deux productions (L2-Strasbourg)}\label{fig:zone3}
\end{center}
\end{figure}
\item
Zone 4 : ce sont les lignes des cas d'égalité et ce sont les moins bien traitées (réussites :  4 sur 10 copies de L2, 13 sur 18 copies en L3 et  2 sur 8 copies de Master). Le fait qu’il n’y ait qu’un nombre fini de solutions (2 solutions) déroute les étudiants. Lorsqu’il y a des réponses avec erreurs, les étudiants arrivent à remplir les colonnes 1 et 2 mais le schéma n'est pas correct. On peut faire l'hypothèse que pour le cas d'égalité, on utilise peu la représentation ensembliste (colonne 1)  dans le secondaire. Le cas d'égalité est un indicateur  de la compréhension de la notion de distance sur la droite réelle.

\begin{figure}[!ht]\label{fig:ex2alg}
\begin{center}
\includegraphics[width=14cm]{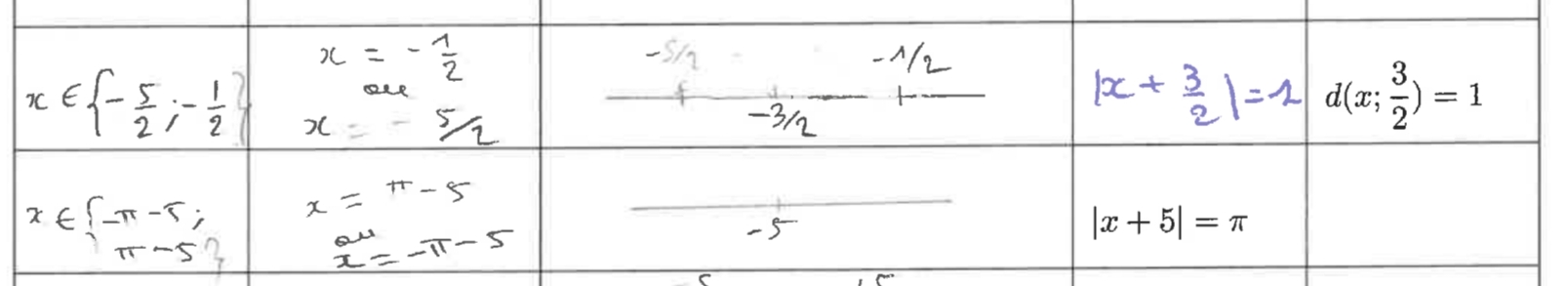}
\caption{Zoom zone 4 : une production (M1-Poitiers)}
\end{center}
\end{figure}
% {\color {purple}
% { Faire réference au cercle- Chevallard transposition didactique  Deuxième édition )}} \cite{CJ}.
\end{itemize}

\subsection{Conclusion pour l'exercice 1}

L'ensemble du tableau est rarement correctement rempli. Les lignes des cas d'égalité sont les moins bien traitées (zone 4) puis viennent les zones 2-1, 2-2 et 3 en raison de la méconnaissance du lien entre valeur absolue et distance. Nous pouvons faire l'hypothèse que les pratiques majoritaires en  calcul ne proposent pas assez souvent l'utilisation de la valeur absolue. On peut imaginer, à la sortie du lycée une prise en charge en L1 d'un travail sur les différents registres mettant en jeu les notions de valeur absolue, distance sur la droite réelle, ensemble de valeurs et  intervalle, pour améliorer la capacité des étudiants à utiliser les changements de registres nécessaires et rendre ces notions opérationnelles.  Dans les exercices portant sur la r\'esolution d'in\'equations, il s'agit de ne pas se contenter d'une r\'esolution alg\'ebrique : \'ecrire l'ensemble des  solutions en utilisant des accolades, ou encore une repr\'esentation sur la droite des r\'eels doit être clairement demandée. 

Nous pouvons  imaginer d'autres exercices en dimension 2 qui sont parfois plus faciles à appréhender, comme des représentations de zones du plan définies à partir  de valeurs absolues ou de distances. 
%\newpage 

\section{Preuve graphique d'une proposition algébrique implicative}
\subsection{Présentation de l'exercice 2 et analyse {\it a priori} }

 Il s'agit ici de démontrer, en utilisant une preuve graphique l'implication suivante  : pour tout r\'eel $a >0$ et pour tout $x \in \R$ 
 
 $$ \left(( \vert x \vert \leq \vert x+a \vert) \mbox { ET } ( \vert x \vert \leq \vert x-a \vert) \right )\Longrightarrow \vert x \vert \leq \frac a 2.$$

Nous avons d'abord testé la production graphique de cette implication lors d'un atelier à destination des enseignants du secondaire. \`A l'issue de cet atelier nous avons choisi de proposer en amont l'exercice de remplissage du tableau précédent, afin de mettre en place le lien entre distance et valeur absolue. 

La preuve graphique est acceptable à condition de préciser les intersections de droites à prendre en compte et les régionnements du plan déterminés par les inégalités.
Des preuves graphiques sont donc possibles en faisant un croquis et en notant clairement les points d'intersection et en complétant  par  des calculs explicatifs : 
\begin{itemize}
\item[a)] en terme de distance  l'hypothèse devient $d(0,x) \leq d(x,-a)$  et  $ d(0,x) \leq  d(x,a) $ et  la représentation  sur la droite de réels de cette situation, peut amener à représenter $-\frac  a 2$  et $\frac a 2$ et à conclure,
\item[b)] à l'aide des représentions graphiques des fonctions $ x \mapsto \vert x\vert$, $x \mapsto \vert x+a \vert $ et $x \mapsto \vert x-a \vert $.
\end{itemize}
Des preuves algébriques sont également possibles :
\begin{itemize}
\item[a)] avec une disjonction de cas, 
\item[b)] en considérant les inégalités équivalentes avec les carrés (tous les termes sont positifs) et en utilisant que $-\frac{a}{2}\leq x\leq \frac{a}{2}$ équivaut à $|x|\leq \frac{a}{2}$.\\
\end{itemize}

\begin{figure}[!ht]
\begin{center}
\includegraphics[width=12cm]{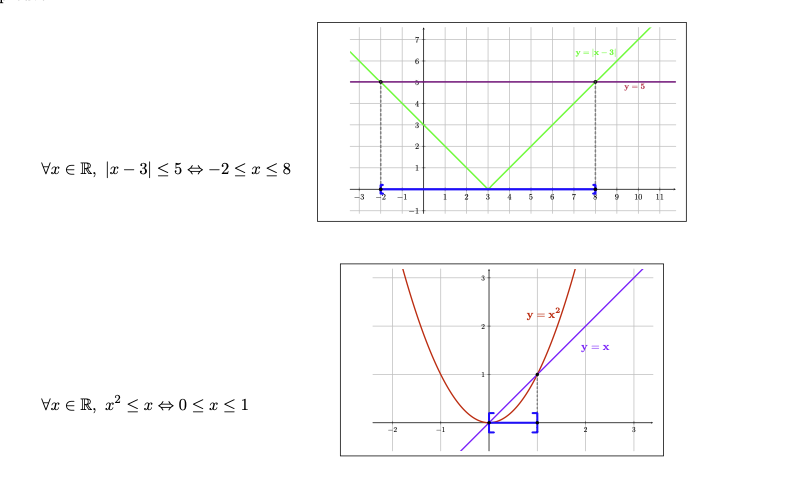} 
\caption{Graphiques ajoutés}\label{fig:exo2}
\end{center}
\end{figure}

Suite à une première expérimentation, nous avons constaté des tentatives de réponses à l'aide des représentations graphiques des trois fonctions, mais  ces représentations étaient portées sur des schémas différents, ce qui empêchait les comparaisons graphiques. Nous avons alors  rajouté les deux graphiques suivants (figure \ref{fig:exo2}). Comme le suggère cette adaptation, le \og ET \fg $\,$  logique de la proposition donnée n'était  pas complètement pris en compte par les étudiants, sans cet ajout.

\vspace{8cm}

Par ailleurs, ces deux graphiques suggèrent le passage d'une proposition  à une inconnue qui peut se résoudre dans $\R$ dans le registre alg\'ebrique, à une représentation dans le plan muni d'un repère, en utilisant un graphique dans le registre fonctionnel. Ce changement de registre est accompagné  par la donnée des exemples précédents.  Nous pouvons faire l'hypothèse que les étudiants utiliseront plutôt une représentation graphique dans $\mathbb R^2$, plutôt que sur la droite réelle en se référant à l'exercice du tableau.

La déduction des tracés de courbes représentatives des fonction $x\mapsto f(x+a)$ et $x\mapsto f(x-a)$ à partir de celle d'une  fonction $f$  (ici $f(x) = \vert x \vert$) semble peu présente dans les pratiques du lycée et peut également être un obstacle dans la mise en oeuvre de la démonstration graphique.

Que ce soit dans le registre graphique ou algébrique, on peut s'attendre à ce que  la présence du paramètre $a$ soit un obstacle à la mise en place d'un raisonnement. Nous avions aussi pour but de voir comment  cet obstacle est contourné, par le choix, par exemple, de valeurs particulières de $a$ pour s'approprier la question.
 
\subsection{Analyse {\it a posteriori}} 

Parmi  les  45  copies (11 de Bordeaux, 7 de Limoges, 8 de Poitiers et 19 de Strasbourg) que nous avons analys\'e 23 ont  d\'ebut\'e une réflexion sur cet exercice : 8 sur les 20 copies de L2, 13 sur les 18 copies de L3 et 2 sur les 8 copies de Master. 
Parmi les étudiants ayant réussi le remplissage du tableau du premier exercice du test, on peut repérer des productions de dessins intéressantes au sujet de l'utilisation du paramètre $a$ mais aussi des échecs. 
Certains essaient  de revenir à la  définition de la valeur absolue, mais la propriété $|x|<a\Leftrightarrow -a<x<a$ n'est pas souvent utilis\'ee (3 fois correctement)  et beaucoup d'étudiants n'utilisent  pas le dessin de la droite réelle comme support de raisonnement.

\begin{figure}[!ht]
\centering
\includegraphics[width=12cm]{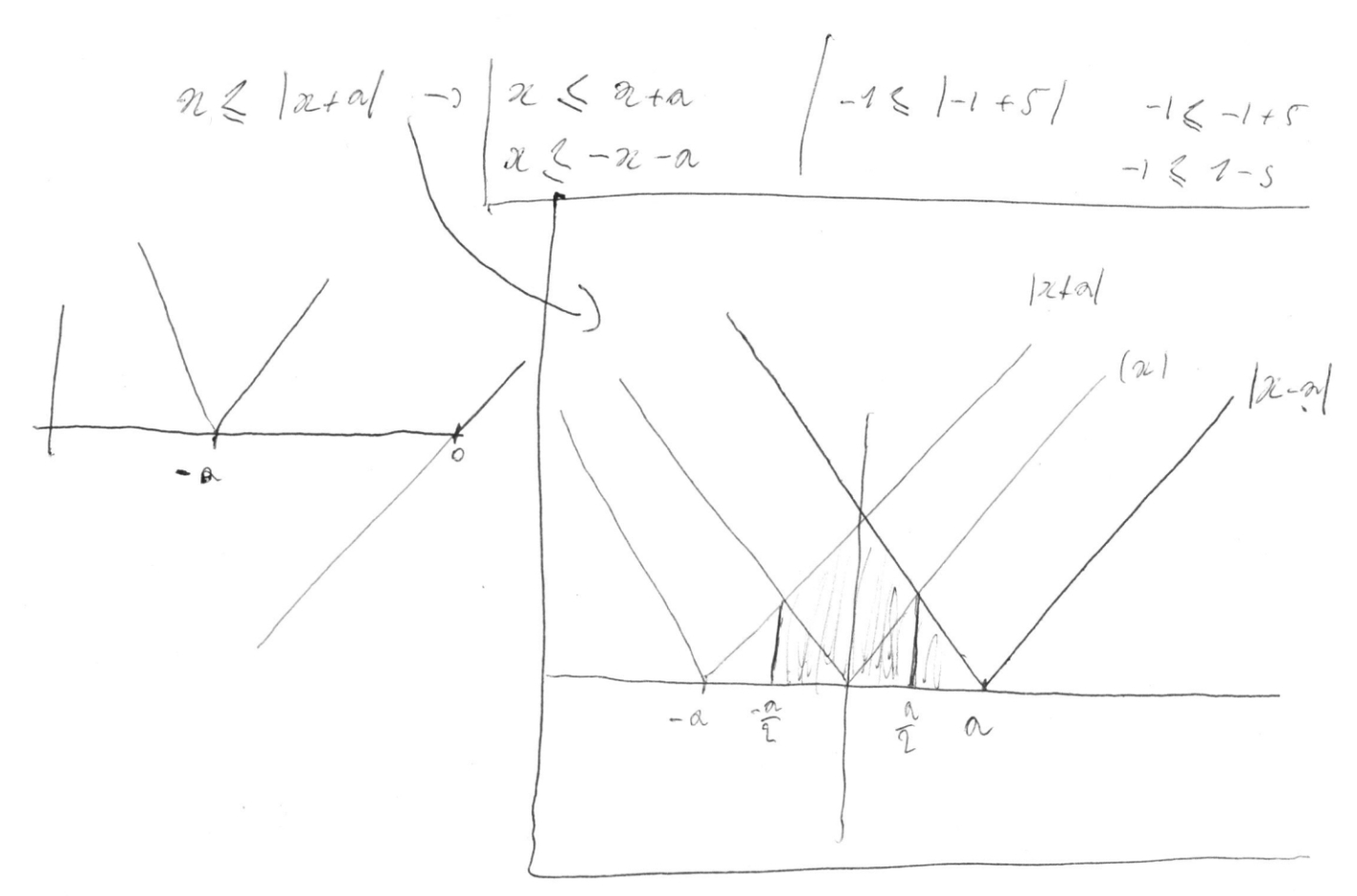} 
\caption{Un exemple de réussite (L2-Strasbourg)}\label{fig:Rparam1}\end{figure}

\begin{figure}[!ht]
\centering
\includegraphics[width=12cm]{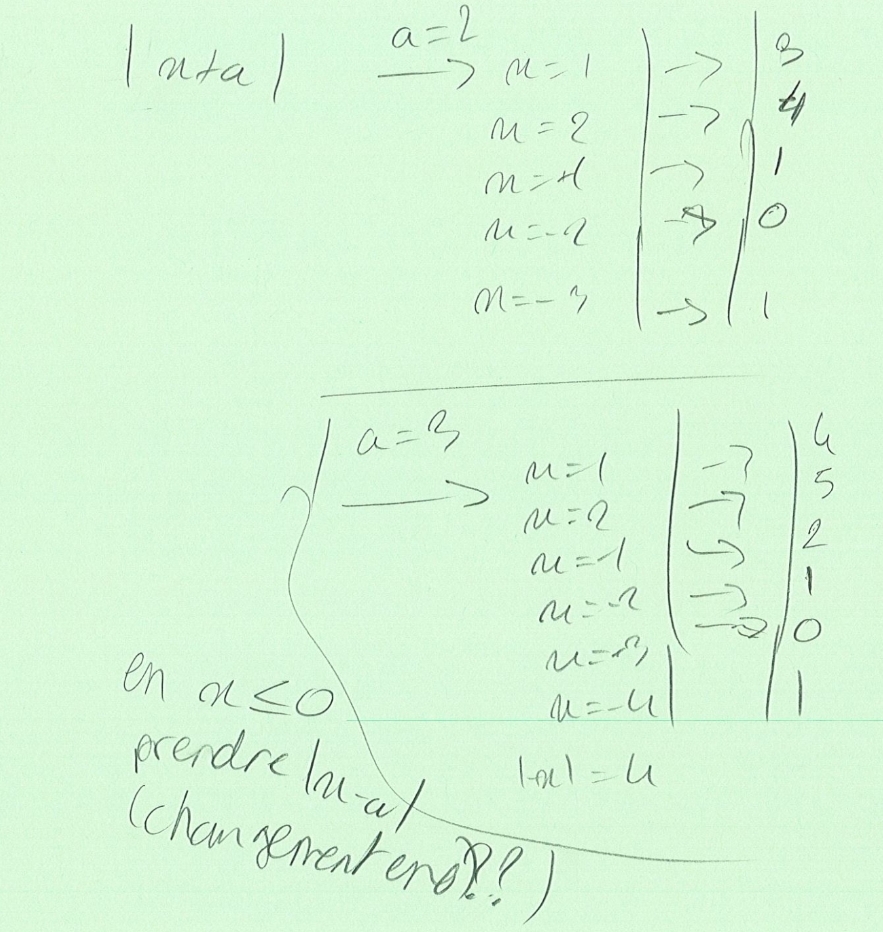} 
\caption{Un exemple d'échec après plusieurs valeurs numériques pour $a$ (L3-Bordeaux)}\label{fig:Eparam2}
\end{figure}

Par ailleurs ci-dessous, deux illustrations avec un brouillon où l'on peut voir un essai de valeur numérique pour $a$ puis le dessin attendu  (figure \ref{fig:Rparam1}) et un autre brouillon où plusieurs valeurs numériques sont prises pour $a$, sans succès (figure \ref{fig:Eparam2}).

\vspace{8cm} 

La déduction des tracés de courbes représentatives des fonctions $x\mapsto f(x+a)$ et $x\mapsto f(x-a)$ à partir de celle de $f$ n'est clairement pas une compétence mobilisable par tous. La  figure \ref{fig:Eparam2}) illustre la difficulté à faire les tracés après plusieurs essais numériques. On peut noter que parmi les 13 copies étudiants de L3, les trois courbes sont représentées sur  9 d'entre elles et sont accompagnées d'une légende correcte mais sans apport significatif pour la résolution. Pour les copies de L2, 1/5 des copies sont de ce type également. Lorsque le graphique proposé est juste, il est rarement accompagné   d'une explication satisfaisante (4 copies sur les 23  ayant débuté une réflexion). 

 Ci dessous deux exemples  de réussite : 
 \begin{figure}[!ht]
\centering
 \includegraphics[width=11cm]{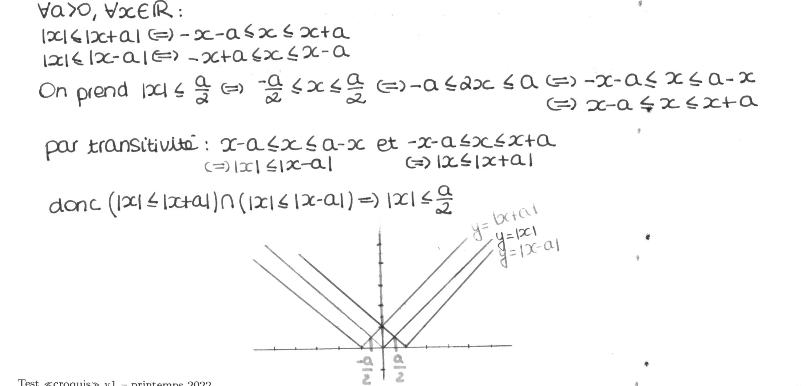} 
 \caption{L2-Strasbourg}\label{fig:exo2strat1}
 \end{figure}
  \begin{figure}[!ht]

\includegraphics[width=12cm]{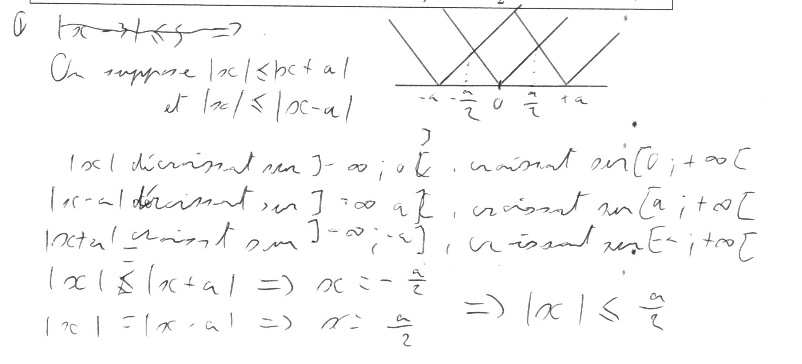} 
\caption{L2-Strasbourg}\label{fig:exo2strat2}

\end{figure}

Dans ces deux productions (figures \ref{fig:exo2strat1} et \ref{fig:exo2strat2}), les trac\'es des trois courbes dans leurs positions respectives sont corrects et l'utilisation du param\`etre ne semble pas être un obstacle. En revanche, le graphique n'est pas clairement utilisé comme support des écritures algébriques qui les accompagnent. En effet, pour que le graphique illustre le résultat à démontrer, il faut identifier les parties du plan correspondant aux deux inégalités de la prémisse et leur conjonction avec le connecteur ET. Cela ne semble pas disponible chez ces étudiants. On peut faire l’hypothèse qu’il y a une double difficulté pour eux : d’une part, lire ces inégalités dans le registre graphique et, d’autre part, les traduire dans le registre algébrique par une disjonction des cas sur la position de la variable $x$ par rapport à $-a$, $0$ et~$a$. 
%sont corrects (les trois courbes sont représentées), sans légende et l'utilisation du paramètre ne semble pas un obstacle pour ces tracés. En revanche, le graphique  n'est pas clairement utilisé comme support des explications algèbriques données (figure \ref{fig:exo2strat2}).
 % 
%Pour utiliser ce graphique il s'agit de traduire graphiquement les deux inégalités de la prémisse. Or, pour les étudiants identifier clairement les parties du graphique à considérer n'est pas immédiat. Ils se retrouvent confrontés à la fois à la  lecture et à la production d'un graphique.
%Une fois les représentations graphiques faites, les justifications sont peut-être évidentes à leur yeux ou bien la preuve graphique ne les guide pas suffisamment pour une preuve algébrique avec par exemple un bon choix pour une disjonction de cas sur la position de $x$ par rapport à $-a$, $0$ et $a$.

\begin{figure}[!ht]
\begin{center}
\includegraphics[width=12cm]{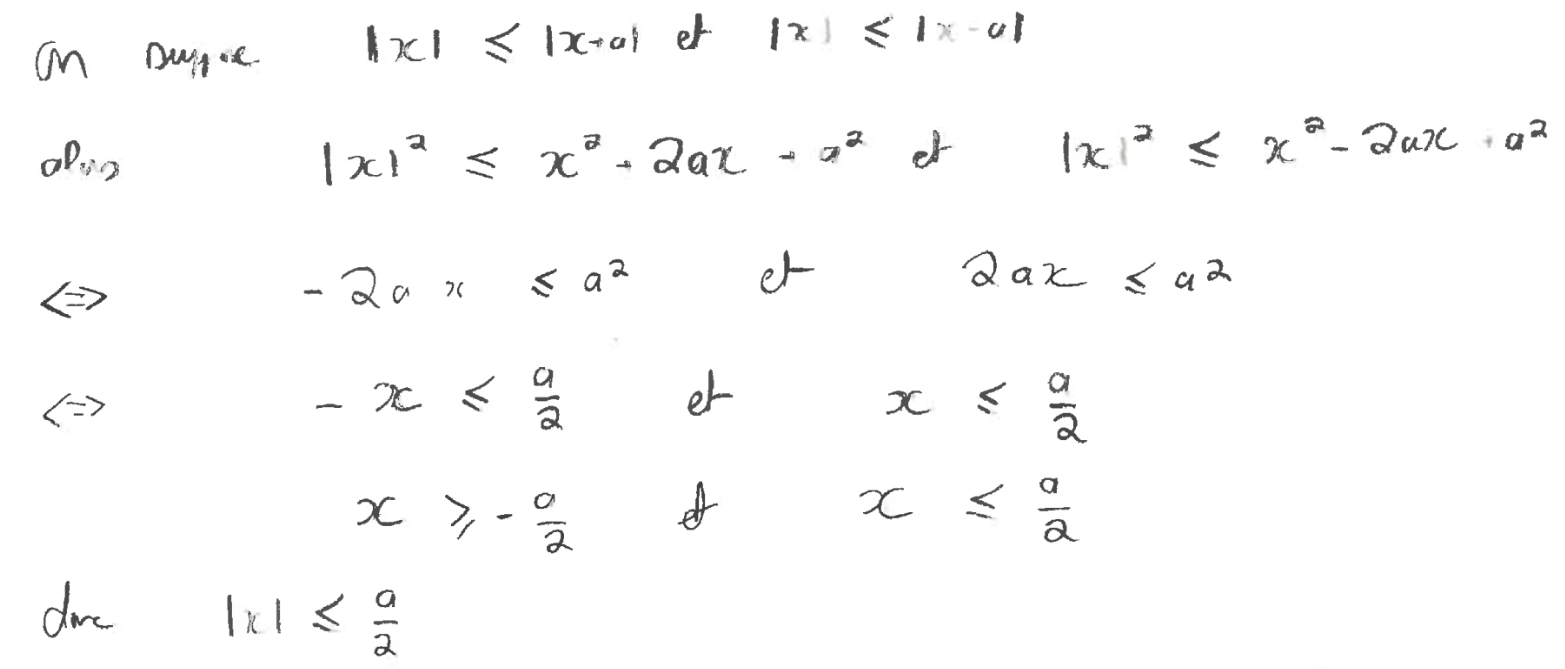}
\caption{Résolution algébrique (M1-Poitiers)}\label{fig:resolalg}
\end{center}
\end{figure}

On peut noter, dans les rédactions proposées, des confusions entre la fonction et son expression et des implications erronées.
Pour le travail algébrique, les inégalités équivalentes au carré sont présentes dans 2 copies, l'une avec une démonstration non aboutie et l'autre réussie (figure \ref{fig:resolalg}). 
\newpage

Nous avons repéré une seule fois l'écriture  de  $\vert x\vert = \sqrt{x^2}$, mais sans aucune utilisation probante. Il 
 Les exemples de lectures graphiques d'inégalités que nous avons ajoutés dans l'exercice (figure \ref{fig:exo2}), ne sont pas accompagnés de la procédure de lecture graphique, comme nous pourrions la réaliser devant des étudiants et l'information apportée est d'autant réduite.

\subsection{Suggestions d'évolution de l'exercice 2}
Si l'on veut suggérer aux étudiants de construire ou lire un graphique sur la droite réelle, on pourrait tenter de renforcer le lien entre les deux premiers exercices du test,  et mettre dans le tableau du premier exercice une ligne avec comme ensemble proposé une demi-droite, pour illustrer l'équivalence (sous l'hypothèse $a\leq b$)~:
$$(x\in ]-\infty, \frac{a+b}{2}])\Leftrightarrow (x\leq \frac{a+b}{2}) \Leftrightarrow (|x-a|\leq |x-b|)\Leftrightarrow d(x;a)\leq d(x;b)\, .$$

On peut aussi ajouter un exercice entre les deux exercices du test dans lequel il serait explicitement demandé de résoudre de deux façons différentes l'inégalité $|x-a|\leq |x-b|$. La traduction de l'inégalité serait alors  demandée en terme d'inégalité de distances puis représentation  sur la droite réelle ; ainsi qu'en terme d'inégalités entre fonctions et représentation dans le plan avec la position relative des courbes représentatives. 
Voici un exemple possible~:\\ 

\noindent \medskip
\fbox{
\begin{minipage}{16.5cm}
On cherche à déterminer les solutions sur $\R$ de l'inéquation
$$
|x+3|\leq |x-5|.
$$
a) Pour cela on introduit les fonctions $f$ et $g$ définies sur $\R$ par $f(x)=|x+3|$ et $g(x)=|x-5|$. Donner dans un même repère les courbes représentatives des deux fonctions et déterminer graphiquement l'ensemble solution en indiquant des tracés utiles pour la résolution.
\\
b) On note $d(a,b)$ la distance entre les réels $a$ et $b$, déterminer les \og ?\fg{} pour que 
$$
|x+3|\leq |x-5| \Leftrightarrow d(x,?)\leq d(x, ?).
$$
Proposer une résolution graphique en travaillant maintenant uniquement sur la droite réelle.
\end{minipage}
}
\medskip

\noindent Nous pourrions également rajouter un exercice de lecture graphique : en donnant un graphique comme dans la figure \ref{fig:exo2}, sans les équivalences écrites sur la gauche, mais en demandant aux étudiants quelles informations ils peuvent tirer d'un tel graphique.

\subsection{Conclusion pour l'exercice 2}
L’énoncé de cet exercice comporte plusieurs éléments à prendre en compte, ceux se situant dans le
registre logique (implication avec prémisse composée) et les autres dans le registre algébrique
(inégalités avec valeurs absolues et paramètre). La consigne étant de donner une preuve graphique
de la proposition, il est donc nécessaire par contrat de changer de registre de représentation, en se
plaçant au choix dans le plan $\R^2$ ou sur la droite réelle. Les exemples donnés dans la deuxième
version du test incitent à choisir le registre $\R^2$. Cependant, l’énoncé est écrit avec une seule variable réelle $x$ et
aucune fonction n'est suggérée, de plus, les étudiants avaient résolu l’exercice 1 concernant des
variables réelles.
\\
Notre objectif est double :
\begin{itemize}
\item vérifier la capacité des étudiants à traduire l’énoncé dans un registre graphique et à décrire les
éléments sur le croquis pour en faire une preuve (par exemple, les coordonnées des points
d’intersection des courbes),
\item montrer aux étudiants comment le registre graphique permet de voir facilement et de
comprendre le résultat à démontrer.
\end{itemize}
Notre deuxième objectif n’a pas vraiment été atteint, car ce problème a révélé des difficultés assez
générales en Licence qui se situent en deçà de la question posée :
\begin{itemize}
\item difficulté à tracer la fonction valeur absolue quand il y a un paramètre (erreur de position de la
courbe par rapport aux axes),
\item difficulté à construire sur un même graphique plusieurs courbes,
\item difficulté de lecture et d’interprétation des inégalités en positions relatives de deux courbes.
\end{itemize}
Même les deux graphiques rajoutés en exemple dans la deuxième version n’ont pas suffi !
Rappelons que l’exercice précédent (exercice 1, cf. § 2), a attesté de la difficulté d'un tiers de ces
étudiants à traduire, sur la droite réelle, la distance entre deux points en valeur absolue (et
inversement). Il semble que même ceux qui ont réussi cette partie de l’exercice n’ont pas
réinvesti cette transcription dans l'exercice 2. Est-ce que, pour eux, tracer la droite réelle n’est pas une
représentation graphique ?
\\
Ces résultats expérimentaux révèlent qu’il est nécessaire de prendre en charge dans les
enseignements de licence les activités de tracés à la main de courbes dans $\R^2$, de lecture et écriture
de données sur la droite réelle, et plus généralement des exercices mettant en relation les registres
algébrique, analytique, graphique et géométrique. Les résultats de nos expérimentations montrent
que ces deux exercices sont tout-à-fait pertinents, à la fois pour construire ces
savoir-faire et pour les évaluer.

\section{Différents registres pour les notions d'injectivité, surjectivité, bijectivité dans le cadre de l'analyse}
\subsection{Présentation des exercices 3, 4 et 5 et analyse {\it a priori}} 

Les exercices 3, 4 et 5 ont pour objectif de tester l’aptitude des étudiants  à changer de registre dans le cadre de l'analyse, sur des notions nouvelles pour des étudiants de  L1\footnote{Première année universitaire},  que sont l'injectivité, la surjectivité et la bijectivité d'une application. 

Les notions d'injectivité, surjectivité, bijectivité sont très importantes dès le début de l'université, en particulier en algèbre linéaire pour étudier les applications linéaires. L'étude des transformations géométriques en tant que telles ayant disparue au lycée, nous avons fait le choix de travailler dans le cadre de l'analyse, plus familier des élèves, dans lequel ils ont déjà rencontré la notion de bijection (théorème de la bijection). En particulier, nous avons fait le choix dans ce test de ne pas travailler avec les \og patates\fg{} représentant des ensembles finis. C'est une représentation usuelle et qui est pourtant souvent utilisée lors de l'introduction de la définition des notions d'injectivité, surjectivité et bijectivité. Dans le cadre de l'analyse, l'illustration de ces notions est assez différente mais il est fondamental que les étudiants les comprennent pourtant aussi dans ce cadre. 
Par ailleurs nous sommes conscients des difficultés liées à la notion de  fonction  en particulier les confusions entre ensemble de départ et d'arrivée. En ce qui concerne le langage mathématique, la difficulté que peut poser les symboles $\exists$ et $\forall$ à un étudiant de ce niveau, peut aussi tout à fait être un obstacle pour qu'il s'exprime clairement dans le registre formel.

En introduction de ces exercices, nous nous sommes inspirés d'un tableau de Rogalski \cite{ROG1} afin d'illustrer ce que nous appelons registres (Nom de la propriété, Langage formalisé, Langage usuel, Équation, Exemple de représentation graphique). Pour les exemples de représentations graphiques nous avons représenté des fonctions discontinues afin de dissocier les notions d'injectivité, surjectivité et bijectivité de celle de continuité.

%%%tableau
\medskip
\noindent
\resizebox{17cm}{!}{
\begin{tabular}{|p{2.8cm}|p{2.45cm}|p{2.6cm}|p{2.8cm}|p{6cm}|}
\hline
\begin{minipage}{2.8cm}\vspace*{0.2cm} Nom de la \\propriété \vspace*{0.2cm}\end{minipage} & \begin{minipage}{2.45cm}\vspace*{0.2cm} Langage \\formalis\'e \vspace*{0.2cm}\end{minipage} & {Langage usuel} & {Equation} & \begin{minipage}{6cm}\vspace*{0.2cm} Exemple de repr\'esentation \\graphique \vspace*{0.2cm}\end{minipage} \\
\hline 
\begin{minipage}{0.98\linewidth}
$f$ est {\bf injective} de $E$ dans $F$
\end{minipage} &
\begin{minipage}{0.98\linewidth}
$\forall x,x'\in E$
si\\ $f(x)=f(x')$ \\alors $x=x'$
\end{minipage}
 &
\begin{minipage}{0.98\linewidth} Tout \'el\'ement\\ $m$ de $F$ a\\ \emph{au plus} un antécédent 
\end{minipage}
&
\begin{minipage}{0.98\linewidth} $\forall m\in F$\\ l'équation $f(x)=m$\\ a \emph{au plus} une solution dans $E$ 
\end{minipage}
&
\begin{minipage}{0.98\linewidth}
\vspace*{0.2cm}
$\forall m\in F$, la droite d'équation $y=m$ coupe $G(f)$ en \emph{au plus} un point.\\
\includegraphics[scale=0.1]{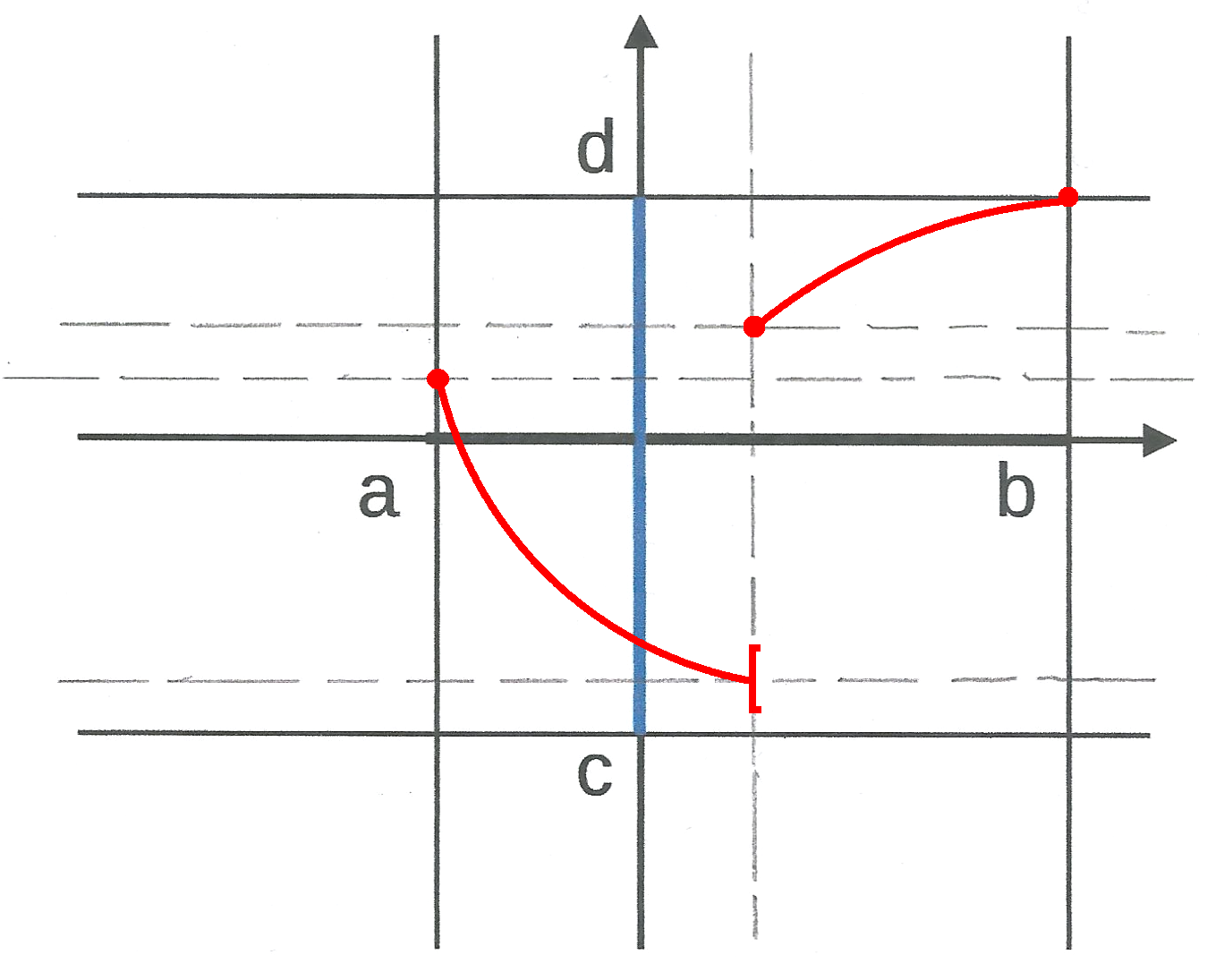}
\\
\footnotesize{avec $E=[a,b]$ et $F=[c,d]$}
\vspace*{0.2cm}
\end{minipage}
\\
\hline 
\begin{minipage}{0.98\linewidth}
$f$ est {\bf surjective} de $E$ dans $F$
\end{minipage} &
\begin{minipage}{0.98\linewidth} $\forall m\in F$\\ $\exists x\in E$ tel que $f(x)=m$
\end{minipage}
 &
\begin{minipage}{0.98\linewidth} Tout \'el\'ement\\  $m$ de $F$ a\\ \emph{au moins} un antécédent 
\end{minipage}
&
\begin{minipage}{0.98\linewidth} $\forall m\in F$\\ l'équation $f(x)=m$\\ a \emph{au moins} une solution dans $E$ 
\end{minipage}
&
\begin{minipage}{0.98\linewidth}
\vspace*{0.2cm}
$\forall m\in F$, la droite d'équation $y=m$ coupe $G(f)$ en \emph{au moins} un point.\\
\includegraphics[scale=0.1]{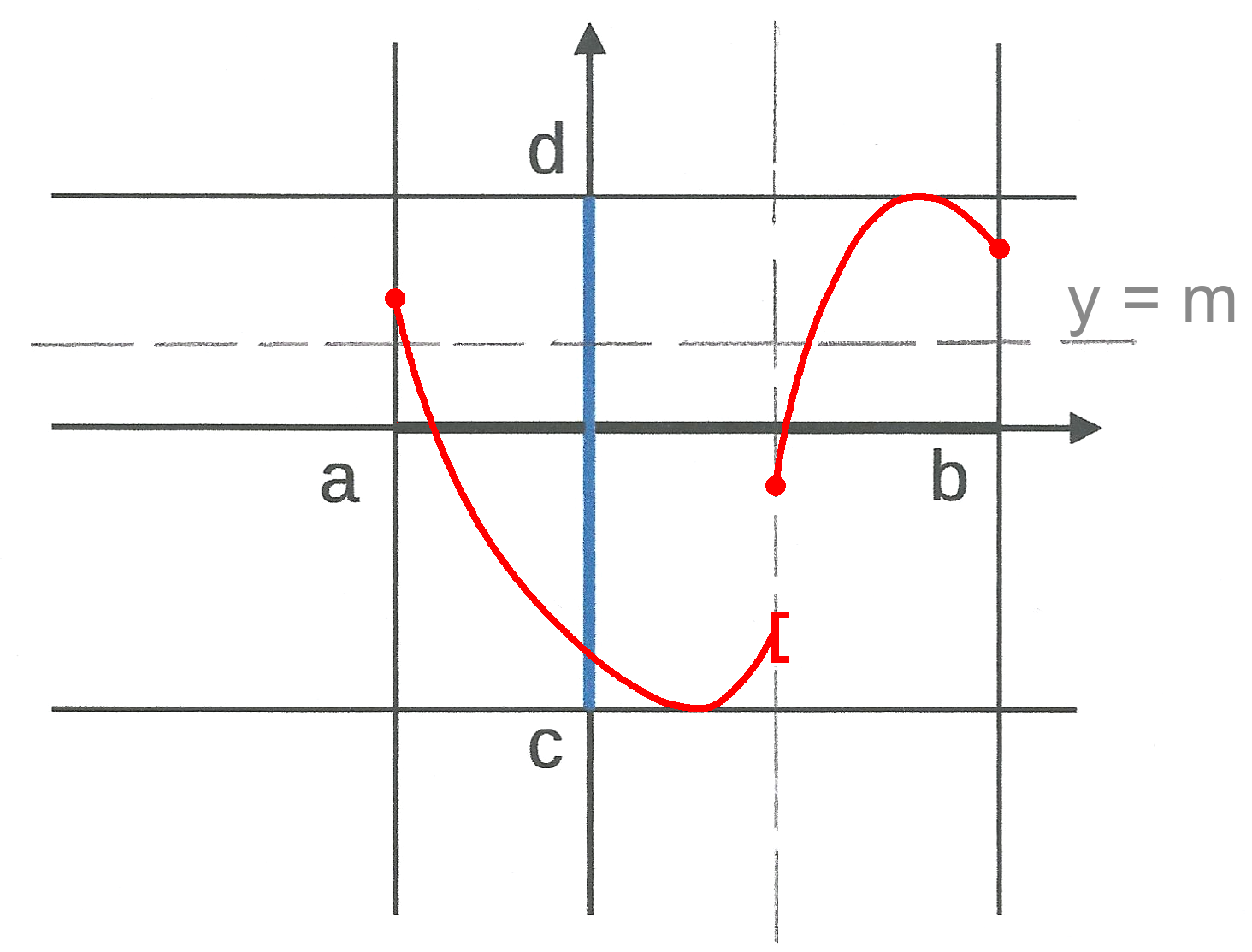}
\\
\footnotesize{avec $E=[a,b]$ et $F=[c,d]$}
\vspace*{0.2cm}
\end{minipage}
\\
\hline 
\begin{minipage}{0.98\linewidth}
$f$ est {\bf bijective} de $E$ dans $F$
\end{minipage} &
\begin{minipage}{0.98\linewidth} $\forall m\in F$\\ $\exists x\in E$ \emph{unique} tel que\\ $f(x)=m$
\end{minipage}
 &
\begin{minipage}{0.98\linewidth} Tout \'el\'ement\\ $m$ de $F$ a \emph{un et un seul} antécédent 
\end{minipage}
&
\begin{minipage}{0.98\linewidth} $\forall m\in F$\\ l'équation $f(x)=m$\\ a \emph{une} solution \emph{et une seule}\\ dans $E$ 
\end{minipage}
&
\begin{minipage}{0.98\linewidth}
\vspace*{0.2cm}
$\forall m\in F$, la droite d'équation $y=m$ coupe $G(f)$ en \emph{un} point \emph{et un seul}.\\
\includegraphics[scale=0.3]{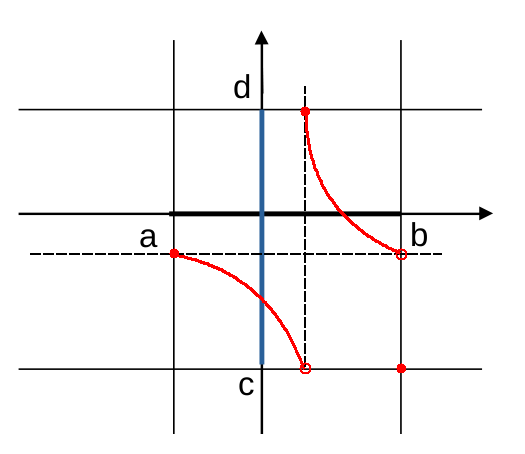}
\\
\footnotesize{avec $E=[a,b]$ et $F=[c,d]$}
\vspace*{0.2cm}
\end{minipage}
\\
\hline 
\end{tabular}\\ 
}

\vspace{1cm} 

Dans l'exercice 3, une même fonction (la fonction carré) est proposée et ne sont modifiés que les ensembles de départ et d'arrivée. Dans chaque cas, il est demandé si les applications sont ou non injective/surjective/bijective en argumentant dans deux registres dont l'un graphique. Le choix a été fait de garder une même fonction bien connue et quatre questions afin qu'une meilleure attention soit portée aux ensembles de départ et d'arrivée.

Pour l'exercice 3, nous avons fait l'hypothèse que le travail à faire le serait d'abord dans le registre graphique avec le tracé des représentations graphiques et une matérialisation des ensembles de départ et d'arrivée. 

Nous avons souhaité voir, en plus de l'illustration graphique, quelle définition ou propriété caractéristique sont privilégiées pour les notions d'injectivité, surjectivité et bijectivité.
Une étude du vocabulaire utilisé est  sans  doute à mener. Ce n'est pas ici notre propos, mais nous sommes conscients que le vocabulaire et son utilisation ont un impact important sur la compréhension de ces notions ce qui rejaillit certainement  sur les représentations graphiques demandées. De même nous ne nous sommes pas concentrés sur les aspects logiques des justifications données.

Dans l'exercice 4, deux représentations graphiques de fonctions ont été proposées avec la matérialisation des ensembles de départ et d'arrivée, en plus de leur mention dans l'énoncé. L'expression analytique des fonctions n'est pas donnée afin d'inciter les étudiants à travailler dans le registre graphique. Le travail donné dans l'encadré ci-dessous consiste à justifier graphiquement si les fonctions sont ou ne sont pas injectives ou surjectives. 
\begin{center}
\fbox{
\begin{minipage}{15cm}
{\small On a représenté ci-dessous deux fonctions de $E=[-1,3]$ dans $F=[-2,3]$. Pour chacune d'elle indiquez si elle est
\begin{itemize}
\item
injective ?
\item
surjective ?
\end{itemize}
Justifiez graphiquement vos r\'eponses.}
\begin{center}
\begin{tikzpicture}[line cap=round,line join=round,>=triangle 45,x=1.5cm,y=1.5cm,scale=0.6]
\begin{axis}[
x=1.5cm,y=1.5cm,
axis lines=middle,
ymajorgrids=true,
xmajorgrids=true,
xmin=-1.5392592592592618,
xmax=3.515555555555562,
ymin=-2.92740740740742,
ymax=3.8518518518518605,
xtick={-2,-1,...,3},
ytick={-2,-1,...,3},]
\clip(-1.5392592592592618,-2.92740740740742) rectangle (3.515555555555562,3.8518518518518605);
\draw[line width=2.5pt,color=blue,smooth,samples=100,domain=-1:3] plot(\x,{2.7*(1/(1+(\x)^2)-0.4)});
\draw [line width=1.5pt,dashed,color=gray] (-1,3)-- (3,3);
\draw [line width=1.5pt,dashed,color=gray] (3,3)-- (3,-2);
\draw [line width=1.5pt,dashed,color=gray] (3,-2)-- (-1,-2);
\draw [line width=1.5pt,dashed,color=gray] (-1,-2)-- (-1,3);
\end{axis}
\end{tikzpicture}
%\end{center}
%
%\vspace{0.2cm}
%
%\begin{center}
\hspace{1.5cm}
\begin{tikzpicture}[line cap=round,line join=round,>=triangle 45,x=1cm,y=1cm,scale=0.5]
%\draw[gray] (-2.7, -4.7) grid (6.7, 6.7);
\draw[thick, black, ->] (-2.7, 0) -- (6.7, 0)
  node[anchor=south west] {$x$};
\draw[thick, black, ->] (0, -4.7) -- (0, 6.7)
  node[anchor=south west] {$y$};
\draw[thick] (-2, 0.2) -- (-2, -0.2) node[anchor=north,scale=1.2] {\($-1$\)};
\draw[thick] (6, 0.2) -- (6, -0.2) node[anchor=north,scale=1.2] {\($3$\)};
\draw[thick] (0.2, -4) --(-0.2, -4) node[anchor=east,scale=1.2] {\($-2$\)};
\draw[thick] (0.2, 6) -- (-0.2, 6)  node[anchor=east,scale=1.2] {\($3$\)};
\draw [line width=1.5pt,dashed,color=gray] (-2,6)-- (6,6);
\draw [line width=1.5pt,dashed,color=gray] (6,6)-- (6,-4);
\draw [line width=1.5pt,dashed,color=gray] (6,-4)-- (-2,-4);
\draw [line width=1.5pt,dashed,color=gray] (-2,-4)-- (-2,6);
\draw [line width=1.5pt,color=blue] (-2,2.6)-- (0.3,1.4);
\draw [line width=1.5pt,color=blue] (0.3,-2)-- (4,4)--(6,3);
\draw[thick, blue, fill=blue] (-2,2.6) circle (1mm);
\draw[thick, blue, fill=white] (0.3,1.4) circle (1mm);
\draw[thick, blue, fill=blue] (0.3,-2) circle (1mm);
\draw[thick, blue, fill=blue] (6,3) circle (1mm);
\end{tikzpicture}
\end{center}
\end{minipage} }
\end{center}

Pour les exercices 3 et 4, nous pensons que l'étude graphique de l'injectivité et de la surjectivité sera faite à l'aide du tracé d'une droite générique d'équation $y=m$ et de la matérialisation de son intersection (ou pas) avec la courbe représentative de la fonction.
Nous attendons que, de plus, les étudiants fassent apparaitre les abscisses des points d'intersection dans le cas non injectif ; qu'ils fassent apparaitre la non-intersection dans le cas non surjectif. En effet, la présence des abscisses sur le schéma est un marqueur du changement de registre, à la différence de deux solutions faites en parallèle. 

Dans l'exercice 4, nous  pensons que l'absence d'expression analytique, même si elle complique un peu le travail, permet de mettre en évidence la compréhension ou non de la traduction graphique.  On peut remarquer que pour répondre correctement à cet exercice il suffit dans chaque cas de donner un contre-exemple. La visualisation graphique d'un contre-exemple est  sans doute plus aisée que la visualisation du  cas générique.  Nous émettons l'hypothèse  que la non-surjectivité sera facilement matérialisée par une droite horizontale  bien choisie qui ne coupe pas le graphe de la fonction.

L'exercice 5 est un exercice qui demande de la créativité. Il est demandé de représenter une fonction qui est injective et non surjective, surjective et non injective, bijective, avec à chaque fois, les mêmes ensembles de départ et d'arrivée qui sont des intervalles.
Le terme de fonction est utilisé ici dans le cadre de l'analyse et nous pensons qu'il sera compris dans le sens usuel d'application.

Une première version du test (copies de Bordeaux et Limoges) demandait des exemples de fonctions affines par morceaux, il y avait eu peu de réussite. Cette consigne a été supprimée dans la version actuelle du test, laissant ainsi la possibilité d'adapter aux bons intervalles les exemples des énoncés. Les étudiants n'ont pas utilisé cette possibilité. 

En plus de demander de la créativité, cet exercice demande aussi de contrôler deux notions en même temps (existence et unicité), nous avons conscience qu'il est plus difficile que les exercices 3 et 4.

\subsection{Analyse {\it a posteriori }}
%
%{\color{red} Rappel de ce qui avait été pointé lors de réunions, il faut maintenant détailler. Les copies qui avaient été regardées de plus près : LeB (très riche) et leK (source L3 Bordeaux).
%
%Axes d'analyse sur les différences de traitement du registre graphique :
%
%- dit mais ne trace rien
%
%- coordonnées représentées mais sans droite (que peut-on en tirer?)
%
%- tracé d'une droite générique
%
%- appui sur le graphique pour raisonner
%
%- la façon dont est exprimé un raisonnement lorsque c'est particulièrement intéressant
%pourra aussi être souligné, mais pas d'étude systématique}
%
%\medskip
%
%{\color{red}Analyse globale pour les copies de Bordeaux et Limoges, à compléter avec Strasbourg et Poitiers en ajoutant du liant}
Rappelons que les analyses portent sur 45 copies (19 copies de L2 Strasbourg, 11 copies de L3 Bordeaux,  7 copies de L3 Limoges et 8 copies de M1 Poitiers).

Parmi les 45 étudiants, 3 n'abordent aucun des trois exercices. Dans les 34 copies où l'exercice 4 est traité, l'exercice 3 l'est aussi. Dans les 29 copies où l'exercice 5 est traité, l'exercice 4 l'est aussi systématiquement.

L'exercice 3 a bien été traité dans le registre graphique, qui était explicitement demandé. Cependant 3 copies sur 42 ne comporte aucun graphique et 2 copies sur 42 comportent uniquement des graphiques faits au brouillon. Ces derniers contenaient des informations et nous les avons donc prises en compte.

\begin{figure}[!ht]
         \begin{center}
         \includegraphics[width=6cm]{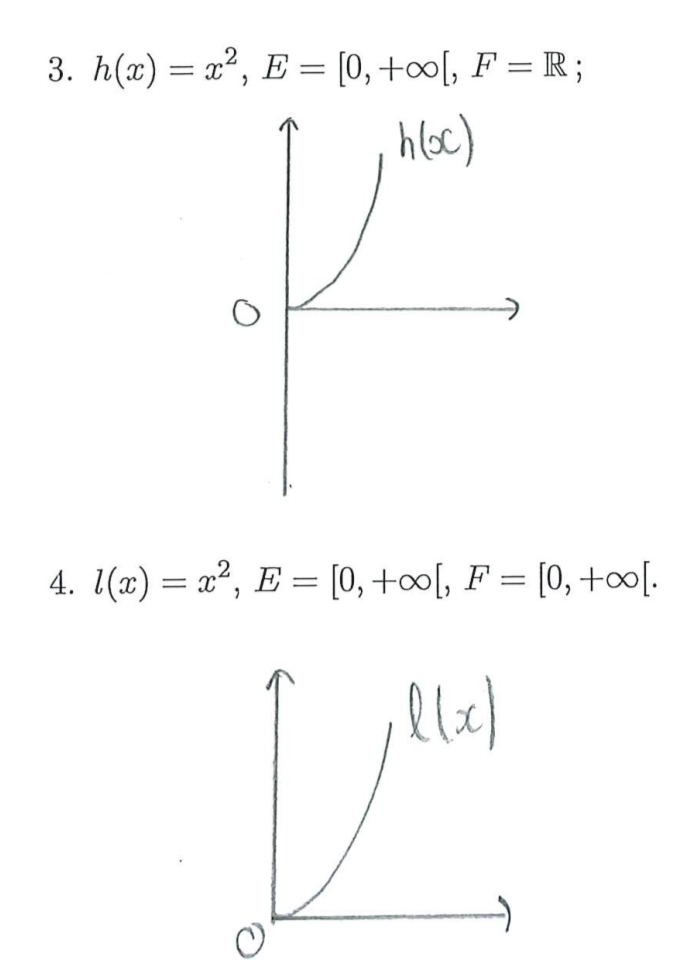}
\caption{Visualisation via les axes de $E$ et $F$ (L3-Limoges)}
\label{fig:axes}
\end{center}
\end{figure}

La matérialisation des ensembles de départ et d'arrivée est faite (33 copies sur 39) et ce le plus souvent via la limitation des axes à des demi-droites lorsque l'un des ensembles est $\R_+$ (figure \ref{fig:axes}). On trouve quelques cas de coloration de la zone du plan correspondant. On trouve aussi une matérialisation partielle avec le tracé d'une demi parabole lorsque l'ensemble de départ $E=\R_+$, qui ne permet pas de matérialiser les différents cas pour l'ensemble d'arrivée.

\newpage
Dans les 39 copies où il y a des graphiques, on note l'absence totale de tracé de  droite horizontale (qu'elle soit particularisée  ou bien générique $y=m$) pour 6 d'entre-elles. On peut alors qualifier la représentation graphique qui a été faite de dessin et non de croquis puisqu'elle n'est pas un support visuel à l'argumentaire.
De plus, le tracé d'une droite générique n'est présent que comme support à l'étude de l'injectivité dans 7 des 33 copies. La droite illustre essentiellement dans les copies la non injectivité (on exhibe une droite horizontale qui coupe 2 fois le graphe) et non surjectivité (on exhibe une droite horizontale qui ne coupe pas le graphe). Il semble qu'il y a chez les étudiants une confusion entre droite particulière à titre de contre exemple et droite générique ($y=m$) qui sert aussi bien de contre exemple "générique" et de support pour justifier l'injectivité ou la surjectivité. 

L'exercice 4 est abordé dans 34 copies. Une droite (particulière ou générique d'équation $y=m$) est tracée sur 30 des 34 copies. On peut supposer qu'il y a une vision générique des croquis (la droite est déplacée verticalement jusqu'à une position permettant de justifier la réponse).

Lorsqu'il y a intersection de cette droite avec la représentation graphique de la fonction étudiée, la matérialisation ou non des abscisses des points d'intersection est faite de façon analogue dans les exercices 3 et 4. Cette matérialisation est présente dans la moitié des copies pour lesquelles nous avons une droite horizontale (16 copies sur 33).

%{\color{red} rem Chantal : S'il y en a besoin il faudra que je regarde mieux dans l'exo 3 le lien entre abscisses placées ou pas et rédaction.}

\begin{figure}[h!]
\centering
\includegraphics[scale=0.7]{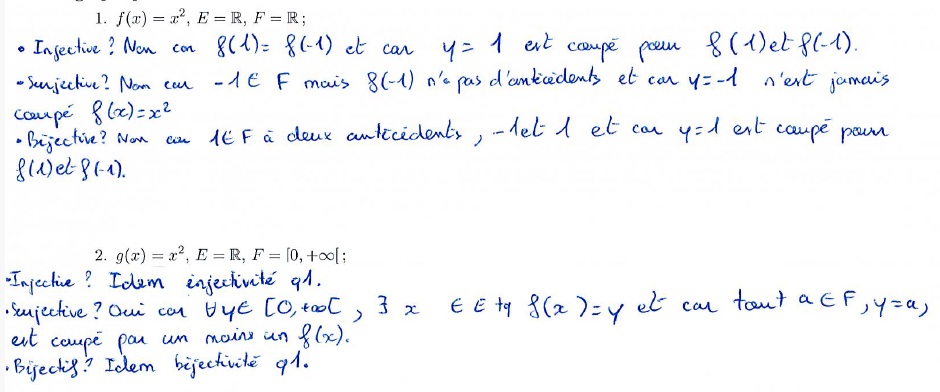}
\caption{Résolution correcte avec un vocabulaire non adapté (L3-Bordeaux)}\label{fig:langage}
\end{figure}

Pour ces deux exercices 3 et 4, bien que tout ne soit pas totalement complété, les étudiants donnent très souvent des bonnes réponses concernant la nature des fonctions même si quelques justifications, notamment dans le registre formel (figure \ref{fig:langage}), sont difficilement acceptables à ce niveau d'études (en L3 notamment).

%\vspace{8cm}

Pour l'exercice 3, nous avons 3 copies sur 42 où il y a la fois des erreurs concernant l'injectivité et la surjectivité, 7 copies sur 42 où il y a des erreurs uniquement sur la surjectivité (dans ce cas le fait que tout point ait une image est confondu avec la notion de surjectivité).
On trouve seulement 2 copies sur 34 erreurs concernant la non-injectivité pour l'exercice 4 et 10 copies sur 34 concernant la non-surjectivité (surtout en L2, 6 copies sur 13 et en M1, 3 copies sur 8).

Les erreurs sur la notion de surjection semble révéler des lacunes sur la notion de fonction :
confusion entre existence et unicité de l'image. Nous savons que de façon générale la différence entre existence (\og au moins un\fg{}) et unicité (\og un et un seul\fg{}) est une difficulté.

On retrouve dans l'exercice 5 que la notion de surjectivité est plus difficile à comprendre que celle d'injectivité, avec des confusions entre existence d'une image (trou dans le graphe), unicité de l'image (condition nécessaire de la relation fonctionnelle) et surjectivité proprement dite.

Dans cet exercice, les trois items ne sont pas tous traités. Il y a échec pour les 3 items dans une copie de L3 et une copie de M1. Nous supposons que la cause en est la prise en compte unique de l'ensemble des images. Dans la copie de L3, on a ainsi à deux reprises des trous dans l'ensemble de départ et le cas d'une représentation qui n'est pas celle d'une fonction (un point avec une infinité d'images). Dans la copie de M1, c'est une non réponse dans l'item {\it injectif et non surjectif} et la non prise en compte de la surjectivité dans les deux autres items. Pour les 27 autres copies~:
\begin{itemize}
\item
l'item {\it injectif et non surjectif} est correctement fait dans seulement 18 copies sur 27 copies (trou dans l'ensemble de départ, dont une représentation avec un seul point isolé),
\item
l'item {\it surjectif et non injectif} est correctement fait dans 21 copies sur 27 copies (un point de l'ensemble d'arrivée n'a pas d'antécédent),
\item
l'item {\it bijectif} est correctement fait dans 24 copies sur 27 copies (en L3 les deux erreurs sont la donnée d'une fonction constante et d'une fonction non injective ; en L2 il s'agit d'un item non renseigné par l'étudiant qui pourtant avait répondu correctement aux deux premiers).
\end{itemize}

\begin{figure}[!ht]
         \centering
 \begin{subfigure}{}
         \includegraphics[width=8cm]{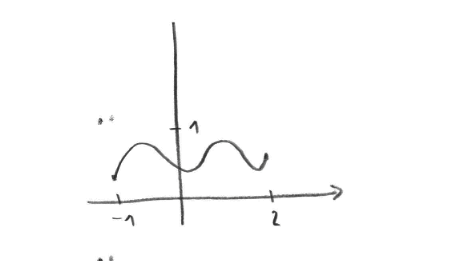}
     \end{subfigure}
     \hfill
     \begin{subfigure}{}
         \includegraphics[width=7cm]{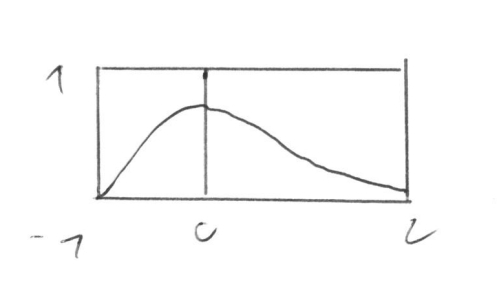}
     \end{subfigure}     
\caption{Exemples pour {\it surjective et non injective} de $[-1,2]$ dans $[0,1]$ (L2-Strasbourg)}
\label{fig:surj pas inj}
\end{figure}

\begin{figure}[!ht]
         \centering
         \includegraphics[width=8cm]{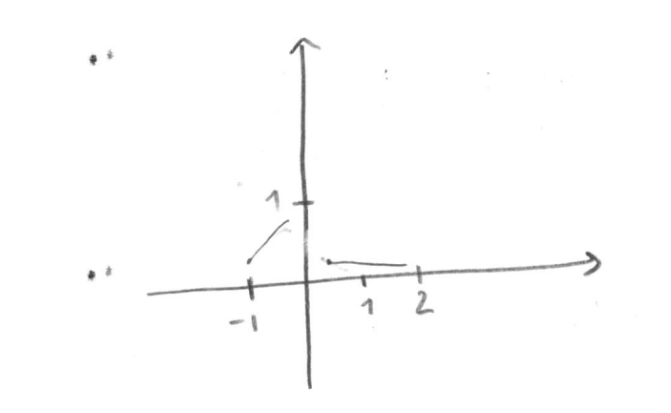}
\caption{Exemple pour {\it injective et non surjective} de $[-1,2]$ dans $[0,1]$ (L2-Strasbourg)}
\label{fig:inj pas surj}
\end{figure}

L'oubli de l'ensemble d'arrivée est une erreur courante. Nous connaissons la difficulté des étudiants à lire les images d'une fonction sur l'axe des ordonnées sans que nous soyons en mesure de l'expliquer. La surjectivité est souvent vue comme \og sans trou dans l'ensemble de départ \fg{} (figures \ref{fig:surj pas inj} et \ref{fig:inj pas surj}). Il est vrai qu'il est demandé dans cet exercice une fonction mais nous pensons que la différence de définition entre application et fonction n'est en général pas prise en compte par les étudiants dans leur pratique dans le cadre de l'analyse. Ce n'est donc pas l'explication pour la présence de ces \og trous \fg{} dans le cas non surjectif, d'autant plus qu'ils ne sont pas présents dans le cas de la non injectivité.

%\vspace{8cm}

%\vspace{8cm}

A cela s'ajoute la confusion entre ensemble d'arrivée et image, ce qui a un impact lorsqu'on demande de produire une fonction non surjective.

Remarquons enfin que l'on peut supposer que pour justifier la non-injectivité, les étudiants ne reviennent pas à la négation de l'injectivité mais ont plutôt la règle en acte : trouver deux points de l'ensemble de départ qui ont la même image.

%{\color{red} Présenter aussi les difficultés pour beaucoup d'étudiant-e-s de nommer/dénommer les objets dont on parle
%(cf. Viviane : Prise en compte de quatre types d'objets et du vocabulaire associé pour équation, fonction, graphe, représentation graphique et leur liens, les trois premiers étant dans le registre numérique, le dernier dans le registre graphique. Par exemple, dans le registre graphique, le mot "point d'intersection" n'apparaît pas ; le
%repérage d'un point par le couple de ses coordonnées non plus ; le mot "antécédent" est utilisé pour parler de solution d'une équation ; il y a un mélange entre vocabulaire géométrique et numérique (exemple - copie LeK, exo 3 : "1 dans $F$ a deux antécédents -
%1 et 1 car $y = 1$ est coupé par $f(1)$ et $f(-1)$" ; et aussi "non car $y = 1,2$ est coupé deux fois
%par $f(x)$"). De nombreux étudiant-e-s étiquètent la représentation graphique avec "$f(c)$".)}

%{\color{red} (Uniquement A partir des copies de Limoges)
%
%- ex 3 : 
%Parfois confusion « image » et « antécédent » dans le vocabulaire : « -1 n’a pas d’image donc pas
%surjectif ».
%Confusion qui relève sans doute des erreurs entre abscisses et ordonnées.
%Sur certaines productions, absence ou mauvais usage des quantificateurs.
%Des difficultés pour prendre la négation d'une proposition en particulier pour une implication.
%On retrouve qu’il est dur de nier injectif/surjectif.}
%
%{\color{red}Peut-être encore des choses à piocher dans le CR du 03/12/21. concernant les copies de Limoges et Bordeaux
%\url{https://plmbox.math.cnrs.fr/f/f0b451dc341242aeab0f/}}

\subsection{Suggestions d'évolution des exercices 3-4-5}

Dans le tableau donné dans l'exercice 3, la colonne \og langage usuel\fg{} peut être modifiée. Une suggestion pour {\it surjectif} serait : \og tout élément de $F$ est atteint\fg{}, pour {\it injectif} serait : \og deux éléments différents de $E$ ont des images différentes \fg{}. De plus la représentation graphique dans le cas {\it bijectif} avec des ronds vides pour les points non atteints prête à confusion car cela peut donner l'impression de ne plus avoir la représentation d'une fonction en ayant deux images pour un même point. Il vaut mieux replacer ces ronds vides par des crochets comme cela a été fait dans les deux cas précédents.
\\
Toujours dans le tableau, afin de montrer que les \og trous \fg{} liés à la discontinuité des fonctions choisies pour illustrer les notions d'injectivité, surjectivité et bijectivité ne sont pas important pour l'illustration graphique de ces notions, on pourrait pour chaque notion présenter deux types de fonctions, l'une discontinue et l'autre continue.

Avant de proposer des exercices tels que les exercices 3-4-5,
il faudrait peut-être commencer avec un exercice uniquement d'observation. Il pourrait prendre la forme suivante~: 4 triplets
$(E,F,f)$ sont représentés, un de chaque type (injectif et pas surjectif, injectif et surjectif, etc ...), et il est demandé d'identifier la nature de chacun en mettant en évidence ce qui permet de décider. 

Les domaines $E$ et $F$, lorsqu'ils
n'étaient pas $\R$ débutaient en 0, de plus, il est bien connu à ce niveau d'étude que la fonction carré est positive sur $\R$. Ces deux points ont peut-être  contribué à ce que les ensembles de départ et d'arrivée ne soient pas matérialisés. On pourrait par exemple proposer l'étude de $x\mapsto x^2-7$ avec $F=\R$ ou $F=[-7,+\infty [$.

\begin{figure}[h!]
\centering
\includegraphics[scale=0.5]{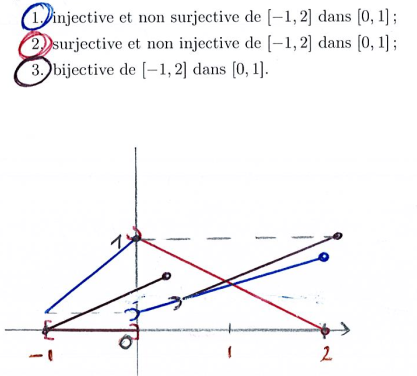}
\caption{Copie L3-Bordeaux}\label{fig:erreur}
\end{figure}

On peut aussi ajouter un exercice avec une ou des erreurs et demander de modifier les tracés, avec un minimum de corrections, afin de les rendre justes. A titre d'exemple, voici une proposition de dessin avec erreurs, pris parmi les productions étudiantes (figure \ref{fig:erreur}), pour laquelle une telle demande pourrait être faite.

\vspace{8cm}

\subsection{Conclusion pour les exercices 3-4-5}

Nous avons constaté dans ces trois exercices qu'il y a plus d'erreurs dans l'identification de la non surjectivité que dans celle de la non injectivité. Les erreurs sont souvent dues à la confusion entre ensemble d'arrivée et ensemble image (droite générique tracée uniquement pour l'étude de l'injectivité dans 7 copies sur 33 cas pour l'exercice 3, 10 copies sur 34 erreurs pour la non surjectivité contre 2 copies sur 34 pour la non injectivité dans l'exercice 4).
On ne travaille probablement pas assez en post-bac sur image/antécédent.

Les suites récurrentes peuvent être une occasion pour travailler ces notions, parler d'intervalles stables et donc de travailler sur un intervalle de départ et l'intervalle d'arrivée correspondant. On peut aussi utiliser les représentations graphiques de suites récurrentes, les valeurs de $u_n$ prenant successivement le rôle d'antécédents de $u_{n+1}$ et d'images de $u_{n-1}$.
Ce travail pourrait être vu en L1 comme une reprise de ce qui peut se faire en terminale avant d'aborder les notions théoriques d'injectivité, surjectivité, bijectivité. 

Dans le post-bac il est assez peu demandé de tracés aux étudiants, on voit dans l'exercice 5 que c'est une lacune importante.
Cet exercice qui demande de la créativité est aussi plus difficile de ce fait, mais les exemples de tracés qui le précédaient auraient pu être adaptés. On peut donc émettre l'hypothèse que les images mentales de l'injectivité et surjectivité dans le cadre des représentations graphiques des fonctions n'est pas en place pour les étudiants en difficulté sur cet exercice. 

On peut aussi souligner que dans l'exercice 3, bien que le registre graphique était explicitement demandé, il n'y a aucun graphique sur 2 des 42 copies et sur 6 autres il y a une absence totale de droite horizontale particularisée ou générique. 

Ces deux points confirment que la production de tracés a besoin d'être travaillée dans l'enseignement supérieur et que les changements de registre faisant intervenir des tracés ont besoin d'être accompagnés.

\section{Conclusion}

Nous pensons que les croquis et les changements de registre qui les accompagnent peuvent aider les élèves ou étudiants à construire de manière opératoire et pérenne les notions mises en jeu dans l'ensemble de notre test. 

Le croquis mis dans la colonne du milieu de l'exercice 1 a pour but d'aider au changement de registres algébrique et géométrique. Les deux exemples de l'exercice 2 sont là pour suggérer le passage du registre algébrique dans $\R$ au registre graphique dans $\R^2$ dans le cadre analytique. Les croquis du tableau de l'exercice 3 sont un moyen de donner une représentation mentale des notions d'injectivité, surjectivité et bijectivité dans le registre fonctionnel.

Nous constatons cependant que les croquis et schémas proposés sous forme achevée (codage et description) ne sont pas opérationnels pour une partie des étudiants testés. Il est donc important d'expliquer en même temps leur construction, en particulier il est compliqué de rendre compte de l'aspect dynamique de ces croquis. 

Cet aspect dynamique est pourtant suggéré par les flèches pour l'exercice 1 (\og placer le point qui sera le milieu puis se déplacer d'une même longueur à gauche et à droite \fg{}), il est absent visuellement pour les croquis de l'exercice 2 (introduction des fonctions, tracé de leur graphe, recherche de point d'intersection et position relative puis lecture de la condition sur $x$ sur l'axe des abscisses) et pour ceux de l'exercice 3 (droite générique d'équation $y=m$ qui se déplace verticalement).

Il nous semble donc indispensable d'inclure régulièrement un travail sur les croquis dans les exercices proposés aux étudiants. Ce n'est qu'ainsi que le croquis deviendra un outil d'aide pour eux. Ainsi les notions que l'enseignant souhaite illustrer avec des croquis seront mieux comprises et ceux-ci seront aussi des supports pour raisonner et résoudre des problèmes.

L'analyse du registre langagier et son lien avec le travail mathématique est une perspective d'étude. Nous faisons l'hypothèse que les difficultés à nommer les objets mathématiques et à en avoir des représentations mentales sont liées.

\newpage

\begin{appendices}

\section{Fiches du test}\label{test}

PSEUDO :

{\sl Communiquer le brouillon en y indiquant aussi le pseudo choisi.}

\bigskip

Ex1 : Dans chaque colonne du tableau ci-dessous est proposée une écriture ou une représentation d'un ensemble de réels. Compléter toutes les cases du tableau afin d'en donner une représentation \'equivalente.

{\sl Indiquer l'ordre de remplissage.}
%Compl\'etez le tableau suivant en traduisant selon les \'ecritures ou repr\'esentations poropos\'ees en colonne.
\begin{center}
\renewcommand{\arraystretch}{3}
\begin{tabular}{|p{2.3cm}|p{3cm}|p{5.8cm}|p{2.3cm}|p{2.3cm}|}

\hline 

\begin{minipage}{2.3cm}ensemble fini\\ de réels ou \\ intervalle(s) \end{minipage} & \begin{minipage}{3cm}\'egalit\'es ou \\ in\'egalit\'es \end{minipage} & \begin{minipage}{7cm}croquis ou sch\'ema de \\ l'ensemble des $x$ considérés \end{minipage} & \begin{minipage}{2.3cm}valeur \\ absolue \end{minipage} & distance  \\ 

\hline 

$x\in]6,15[$ & &  &   & \\ 

\hline 

 & $-\sqrt{2}\leq x\leq \sqrt{2}$ &  &   & \\ 
 
 \hline

 &  & 
\begin{tikzpicture}[line cap=round,x=2cm,y=2cm, scale=0.71]
%\clip(0,-3) rectangle (6,1);
\draw[dashed,color=gray] (1,0)--(5,0);
\draw[line width=1.5,color=blue] ((1.5,0)--(4.5,-0);
\draw[line width=1.5,color=blue] (1.45,0.1)--(1.5,0.1)--(1.5,-0.1) node[below] {$-11$}--(1.45,-0.1);
\draw[line width=1.5,color=blue] ((4.55,0.1) --(4.5,0.1)--(4.5,-0.1) node[below] {$0$}--(4.55,-0.1);
\draw[line width=1,color=gray] (3,0.1)--(3,-0.1) node[below] {\color{blue}$-\frac{11}{2}$};

\draw[color=red,-triangle 60] ((3,.2) to[out =150, in=60] node[midway,above] {$-\frac{11}{2}$} (1.5,.2);
\draw[color=red,-triangle 60] ((3,.2) to[out =30, in=120] node[midway,above] {$+\frac{11}{2}$} (4.5,.2);

\end{tikzpicture}
  &  &  \\ 

\hline 

& & & & $d(x;\frac{3}{2})=1$\\

\hline

& & & $|x+5|=\pi$ & \\

 \hline 

&  &  & & $d(x;-4)\leq 5$ \\ 

\hline

 &  &  & $|x-3|\leq \dfrac{1}{2}$  & \\ 
 
\hline

& & 
\begin{tikzpicture}[line cap=round,x=2cm,y=2cm, scale=0.71]
%\clip(0,-3) rectangle (6,1);
\draw[dashed,color=gray] (1,0)--(5,0);
\draw[line width=1.5,color=blue] ((1,0)--(2,-0);
\draw[line width=1.5,color=blue] ((4,0)--(5,-0);
\draw[line width=1.5,color=blue] (1.95,0.1)--(2,0.1)--(2,-0.1) node[below] {$-\frac{5}{3}$}--(1.95,-0.1);
\draw[line width=1.5,color=blue] ((4.05,0.1) --(4,0.1)--(4,-0.1) node[below] {$\frac{7}{3}$}--(4.05,-0.1);
\draw[line width=1,color=gray] (3,0.1)--(3,-0.1) node[below] {\color{blue}$\frac{1}{3}$};

\draw[color=red,-triangle 60] ((3,.2) to[out =150, in=60] node[midway,above] {$-2$} (2,.2);
\draw[color=red,-triangle 60] ((3,.2) to[out =30, in=120] node[midway,above] {$+2$} (4,.2);

\end{tikzpicture}
 & & \\

\hline 
\end{tabular} 
\end{center}

\newpage

Ex2 : Voici deux exemples de graphiques illustrant la proposition à leur gauche et pouvant servir de preuve.

%{\color{red} Pb avec le dessin des points d'intersection, ne comprend pas ... mis en commentaire dans le .tex}

\begin{center}

\begin{tabular}{cc}
\begin{minipage}{6cm}

$\forall x\in \R, \ |x-3|\leq 5\Leftrightarrow -2\leq x\leq 8$

\vspace*{2cm}

\end{minipage}
&
\fbox{
\definecolor{qqqqff}{rgb}{0,0,1}
\definecolor{uuuuuu}{rgb}{0.26666666666666666,0.26666666666666666,0.26666666666666666}
\definecolor{yqqqyq}{rgb}{0.5019607843137255,0,0.5019607843137255}
\definecolor{qqffqq}{rgb}{0,1,0}
\begin{tikzpicture}[line cap=round,line join=round,>=triangle 45,x=1cm,y=1cm,scale=0.5]
\begin{axis}[
x=1cm,y=1cm,
axis lines=middle,
ymajorgrids=true,
xmajorgrids=true,
xmin=-3.426666666666664,
xmax=11.746666666666659,
ymin=-1.1288888888888857,
ymax=7.564444444444428,
xtick={-3,-2,...,11},
ytick={-1,0,...,7},]
\clip(-3.426666666666664,-1.1288888888888857) rectangle (11.746666666666659,7.564444444444428);
\draw[line width=2pt,color=qqffqq,smooth,samples=100,domain=-3.426666666666664:11.746666666666659] plot(\x,{abs((\x)-3)});
\draw [line width=2pt,color=yqqqyq,domain=-3.426666666666664:11.746666666666659] plot(\x,{(5-0*\x)/1});
\draw [line width=2pt,dashed,color=gray] (-2,5)-- (-2,0);
\draw [line width=2pt,dashed,color=gray] (8,5)-- (8,0);
\draw [line width=3pt,color=qqqqff] (-1.9,0.2)--(-2,0.2)--(-2,-0.2)--(-1.9,-0.2) (-2,0)-- (8,0) (7.9,0.2)--(8,0.2)--(8,-0.2)--(7.9,-0.2);
\draw (9.213333333333328,5.057777777777766) node[anchor=north west]{$\mathbf{\color{purple} y=5}$};
\draw [color=qqffqq](6.946666666666662,6.711111111111096) node[anchor=north west] {$\mathbf{\color{green} y=|x-3|}$};
%\begin{scriptsize}
%\draw [fill=uuuuuu] (-2,5) circle (2pt);  %%%pb incompréhensible
%\draw [fill=uuuuuu] (8,5) circle (2pt);
%\draw [fill=black] (-2,0) circle (2pt);
%\draw [fill=black] (8,0) circle (2pt);
%\end{scriptsize}
\end{axis}
\end{tikzpicture}
}
\\
\begin{minipage}{6cm}
$\forall x\in \R, \ x^2\leq x\Leftrightarrow 0\leq x\leq 1$

\vspace*{2cm}
\end{minipage}
&
\fbox{
\definecolor{qqqqff}{rgb}{0,0,1}
\definecolor{uuuuuu}{rgb}{0.26666666666666666,0.26666666666666666,0.26666666666666666}
\definecolor{xfqqff}{rgb}{0.4980392156862745,0,1}
\definecolor{ccqqqq}{rgb}{0.8,0,0}
\begin{tikzpicture}[line cap=round,line join=round,>=triangle 45,x=2.5cm,y=2.5cm,scale=0.45]
\begin{axis}[
x=2.5cm,y=2.5cm,
axis lines=middle,
ymajorgrids=true,
xmajorgrids=true,
xmin=-2.4049382716049394,
xmax=3.396543209876541,
ymin=-0.5837037037037022,
ymax=3.185185185185181,
xtick={-2,-1,...,3},
ytick={-1,0,...,3},]
\clip(-2.4049382716049394,-0.5837037037037022) rectangle (3.396543209876541,3.185185185185181);
\draw[line width=2pt,color=ccqqqq,smooth,samples=100,domain=-2.4049382716049394:3.396543209876541] plot(\x,{(\x)^(2)});
\draw [line width=2pt,color=xfqqff,domain=-2.4049382716049394:3.396543209876541] plot(\x,{(-0--1*\x)/1});
\draw [line width=2pt,dashed ,color=gray] (1,1)-- (1,0);
\draw (1.826172839506171,1.7) node[anchor=north west,color=xfqqff,scale=1.5] {$ \mathbf{y=x}$};
\draw (0.7,2.462222222222219) node[anchor=north west,color=ccqqqq,scale=1.5] {$\mathbf{ y=x^2}$};
\draw [line width=3pt,color=qqqqff] (0.1,0.2)--(0,0.2)--(0,-0.2)--(0.1,-0.2) (0,0)-- (1,0) (0.9,0.2)--(1,0.2)--(1,-0.2)--(0.9,-0.2);
%\begin{scriptsize}
%\draw [fill=black] (1,1) circle (2pt);    %%%pb incompréhensible
%\draw [fill=black] (1,0) circle (2pt);
%\draw [fill=black] (0,0) circle (2pt);
%\end{scriptsize}
\end{axis}
\end{tikzpicture}
}
\end{tabular}
\end{center}

\medskip

Donnez une démonstration de la proposition suivante et faites un graphique illustrant cette démonstration ou pouvant la remplacer.

%\medskip

\fbox{
\begin{minipage}{15.7cm}
Pour tout r\'eel $a>0$, pour tout $x\in\R$, 
$$ 
\Big((|x|\leq |x+a|)\mbox{ ET }(|x|\leq |x-a|)\Big)\Rightarrow |x|\leq \frac{a}{2}.
$$
\end{minipage}
}

\vspace*{7cm}

\newpage

Ex3 : Le tableau ci-dessous présente les notions d'injection, surjection et bijection dans différents registres (Nom de la propriété, Formel, Langage usuel, Equation, Exemple de repr\'esentation graphique).

\medskip
\noindent
\resizebox{17cm}{!}{
\begin{tabular}{|p{2.8cm}|p{2.45cm}|p{2.6cm}|p{2.8cm}|p{6cm}|}
\hline
\begin{minipage}{2.8cm}\vspace*{0.2cm} Nom de la \\propriété \vspace*{0.2cm}\end{minipage} & \begin{minipage}{2.45cm}\vspace*{0.2cm} Langage \\formalis\'e \vspace*{0.2cm}\end{minipage} & {Langage usuel} & {Equation} & \begin{minipage}{6cm}\vspace*{0.2cm} Exemple de repr\'esentation \\graphique \vspace*{0.2cm}\end{minipage} \\
\hline 
\begin{minipage}{0.98\linewidth}
$f$ est {\bf injective} de $E$ dans $F$
\end{minipage} &
\begin{minipage}{0.98\linewidth}
$\forall x,x'\in E$
si\\ $f(x)=f(x')$ \\alors $x=x'$
\end{minipage}
 &
\begin{minipage}{0.98\linewidth} Tout \'el\'ement\\ $m$ de $F$ a\\ \emph{au plus} un antécédent 
\end{minipage}
&
\begin{minipage}{0.98\linewidth} $\forall m\in F$\\ l'équation $f(x)=m$\\ a \emph{au plus} une solution dans $E$ 
\end{minipage}
&
\begin{minipage}{0.98\linewidth}
\vspace*{0.2cm}
$\forall m\in F$, la droite d'équation $y=m$ coupe $G(f)$ en \emph{au plus} un point.\\
\includegraphics[scale=0.1]{injectivite_copie2.png}
\\
\footnotesize{avec $E=[a,b]$ et $F=[c,d]$}
\vspace*{0.2cm}
\end{minipage}
\\
\hline 
\begin{minipage}{0.98\linewidth}
$f$ est {\bf surjective} de $E$ dans $F$
\end{minipage} &
\begin{minipage}{0.98\linewidth} $\forall m\in F$\\ $\exists x\in E$ tel que $f(x)=m$
\end{minipage}
 &
\begin{minipage}{0.98\linewidth} Tout \'el\'ement\\  $m$ de $F$ a\\ \emph{au moins} un antécédent 
\end{minipage}
&
\begin{minipage}{0.98\linewidth} $\forall m\in F$\\ l'équation $f(x)=m$\\ a \emph{au moins} une solution dans $E$ 
\end{minipage}
&
\begin{minipage}{0.98\linewidth}
\vspace*{0.2cm}
$\forall m\in F$, la droite d'équation $y=m$ coupe $G(f)$ en \emph{au moins} un point.\\
\includegraphics[scale=0.1]{surjectivite_copie2.png}
\\
\footnotesize{avec $E=[a,b]$ et $F=[c,d]$}
\vspace*{0.2cm}
\end{minipage}
\\
\hline 
\begin{minipage}{0.98\linewidth}
$f$ est {\bf bijective} de $E$ dans $F$
\end{minipage} &
\begin{minipage}{0.98\linewidth} $\forall m\in F$\\ $\exists x\in E$ \emph{unique} tel que\\ $f(x)=m$
\end{minipage}
 &
\begin{minipage}{0.98\linewidth} Tout \'el\'ement\\ $m$ de $F$ a \emph{un et un seul} antécédent 
\end{minipage}
&
\begin{minipage}{0.98\linewidth} $\forall m\in F$\\ l'équation $f(x)=m$\\ a \emph{une} solution \emph{et une seule}\\ dans $E$ 
\end{minipage}
&
\begin{minipage}{0.98\linewidth}
\vspace*{0.2cm}
$\forall m\in F$, la droite d'équation $y=m$ coupe $G(f)$ en \emph{un} point \emph{et un seul}.\\
\includegraphics[scale=0.3]{bijection.png}
\\
\footnotesize{avec $E=[a,b]$ et $F=[c,d]$}
\vspace*{0.2cm}
\end{minipage}
\\
\hline 
\end{tabular}
}

\newpage

En vous inspirant de ce qui est proposé dans le tableau précédent, justifiez si les applications ci-dessous sont ou non injective/surjective/bijective {\bf en argumentant dans deux registres dont l'un graphique}.
\begin{enumerate}
\item
$f(x)=x^2$, $E=\R$, $F=\R$ ;

\vspace*{5cm}

\item
$g(x)=x^2$, $E=\R$, $F=[0,+\infty[$ ;

\vspace*{5cm}

\item
$h(x)=x^2$, $E=[0,+\infty[$, $F=\R$ ;

\vspace*{5cm}

\item
$l(x)=x^2$, $E=[0,+\infty[$, $F=[0,+\infty[$.

\end{enumerate}

\newpage

Ex4 : On a représenté ci-dessous deux fonctions de $E=[-1,3]$ dans $F=[-2,3]$. Pour chacune d'elle indiquez si elle est
\begin{itemize}
\item
injective ?
\item
surjective ?
\end{itemize}
Justifiez graphiquement vos r\'eponses.

\begin{center}
\begin{tikzpicture}[line cap=round,line join=round,>=triangle 45,x=1.5cm,y=1.5cm]
\begin{axis}[
x=1.5cm,y=1.5cm,
axis lines=middle,
ymajorgrids=true,
xmajorgrids=true,
xmin=-1.5392592592592618,
xmax=3.515555555555562,
ymin=-2.92740740740742,
ymax=3.8518518518518605,
xtick={-2,-1,...,3},
ytick={-2,-1,...,3},]
\clip(-1.5392592592592618,-2.92740740740742) rectangle (3.515555555555562,3.8518518518518605);
\draw[line width=2.5pt,color=blue,smooth,samples=100,domain=-1:3] plot(\x,{2.7*(1/(1+(\x)^2)-0.4)});
\draw [line width=1.5pt,dashed,color=gray] (-1,3)-- (3,3);
\draw [line width=1.5pt,dashed,color=gray] (3,3)-- (3,-2);
\draw [line width=1.5pt,dashed,color=gray] (3,-2)-- (-1,-2);
\draw [line width=1.5pt,dashed,color=gray] (-1,-2)-- (-1,3);
\end{axis}
\end{tikzpicture}
\end{center}

\vspace{1cm}

\begin{center}
\begin{tikzpicture}[line cap=round,line join=round,>=triangle 45,x=1cm,y=1cm,scale=0.7]
%\draw[gray] (-2.7, -4.7) grid (6.7, 6.7);
\draw[thick, black, ->] (-2.7, 0) -- (6.7, 0)
  node[anchor=south west] {$x$};
\draw[thick, black, ->] (0, -4.7) -- (0, 6.7)
  node[anchor=south west] {$y$};
\draw[thick] (-2, 0.2) -- (-2, -0.2) node[anchor=north,scale=1.2] {\($-1$\)};
\draw[thick] (6, 0.2) -- (6, -0.2) node[anchor=north,scale=1.2] {\($3$\)};
\draw[thick] (0.2, -4) --(-0.2, -4) node[anchor=east,scale=1.2] {\($-2$\)};
\draw[thick] (0.2, 6) -- (-0.2, 6)  node[anchor=east,scale=1.2] {\($3$\)};
\draw [line width=1.5pt,dashed,color=gray] (-2,6)-- (6,6);
\draw [line width=1.5pt,dashed,color=gray] (6,6)-- (6,-4);
\draw [line width=1.5pt,dashed,color=gray] (6,-4)-- (-2,-4);
\draw [line width=1.5pt,dashed,color=gray] (-2,-4)-- (-2,6);
\draw [line width=1.5pt,color=blue] (-2,2.6)-- (0.3,1.4);
\draw [line width=1.5pt,color=blue] (0.3,-2)-- (4,4)--(6,3);
\draw[thick, blue, fill=blue] (-2,2.6) circle (1mm);
\draw[thick, blue, fill=white] (0.3,1.4) circle (1mm);
\draw[thick, blue, fill=blue] (0.3,-2) circle (1mm);
\draw[thick, blue, fill=blue] (6,3) circle (1mm);
\end{tikzpicture}
\end{center}

\newpage

Ex5 : Pour chacune des conditions ci-dessous, repr\'esentez graphiquement une fonction qui la satisfait~:
\begin{enumerate}
\item
injective et non surjective de $[-1,2]$ dans $[0,1]$ ;
\vspace*{7cm}
\item
surjective et non injective de $[-1,2]$ dans $[0,1]$ ;
\vspace*{7cm}
\item
bijective de $[-1,2]$ dans $[0,1]$.
\vspace*{7cm}
\end{enumerate}

\newpage

%Une fonction affine est une fonction dont la repr\'esentation graphique est une droite. Une fonction {\bf affine par morceaux} est une fonction définie sur une réunion d'intervalles et qui est la restriction d'une fonction affine sur chacun de ces intervalles. Voici un exemple de fonction affine par morceaux de $[a,b]$ dans $[c,d]$.

%\section{Travaux d'étudiants}
%
%extraits de travaux d'étudiants qu'on n'aurait pas mis dans le corps du texte et qui ont un intérêt.

\end{appendices}

% Les utilisateurs de BibTeX doivent utiliser
\bibliographystyle{apacite}
\bibliography{CroquisCIIU_ref}
% name your BibTeX data base

\authoraddresses{
Denise Grenier\\
Université Grenoble Alpes, Institut Fourier - IREM\\
100 rue des Maths\\
38 610 Gières\\
\email denise.grenier@univ-grenoble-alpes.fr

Chantal Menini\\
Université de Bordeaux, UF MI - IREM\\
351 cours de la libération\\
33 405 Talence Cedex\\
\email chantal.menini@u-bordeaux.fr

Pascale Sénéchaud\\
Université de Limoges, FST - IREM\\
123 avenue Albert Thomas\\
87 060 Limoges Cedex\\
\email pascale.senechaud@unilim.fr

Fabrice Vandebrouck\\
Université Paris Cité, UFR de Mathématiques - IREM\\
CC 7018\\
75 205 Paris Cedex 13\\
\email vandebro@univ-paris-diderot.fr
}

\end{document}